\documentclass[a4paper]{article}
\usepackage{RR}
\usepackage{hyperref}

\usepackage[english]{babel}
\usepackage[latin1]{inputenc}
\usepackage[T1]{fontenc}
\usepackage{float}
\usepackage{amsmath}
\usepackage{amsfonts}
\usepackage{amssymb}	
\usepackage{amsthm}
\usepackage{stmaryrd}
\usepackage{url}
\usepackage{epsfig}
\usepackage{natbib}
\usepackage{color}
\def\cov{\mbox{cov}}
\def\var{\mbox{Var}}

\def\pen{\mbox{pen}}
\psfigdriver{dvips}
\usepackage{mathrsfs}

\usepackage{geometry}
\newfloat{Figure}{H}{lof}
\newfloat{Table}{Hptb}{lot} 

\newtheorem{thrm}{Theorem}[section]
\newtheorem{prte}[thrm]{Proposition}
\newtheorem{lemma}[thrm]{Lemma}
\newtheorem{cor}[thrm]{Corollary}
\newtheorem{defi}[thrm]{Definition}

\newtheorem{ex}[thrm]{Example}

\RRdate{Janvier 2009}
\RRversion{2}
\RRdater{Septembre  2009}

\RRauthor{Nicolas Verzelen
\thanks{Laboratoire de Math\'ematiques UMR 8628, Universit\'e Paris-Sud, 91405 Osay}
\thanks{INRIA Saclay, Projet SELECT, Universit\'e Paris-Sud, 91405 Osay}}
\authorhead{Verzelen}

\RRtitle{Estimation adaptative de champs gaussiens stationnaires}
\RRetitle{Adaptive estimation of stationary Gaussian fields}
\titlehead{Estimation of stationary Gaussian fields}
\RRresume{Nous étudions l'estimation non-paramétrique d'un champ gaussien stationnaire $X$ observé sur un réseau régulier. Dans le cadre des séries temporelles, certaines procédures comme AIC réalisent une sélection de modèle optimale parmi les modèles autorégressifs. Cependant, il n'existe aucun résultat analogue d'adaptation pour des champs spatiaux. En considérant des collections de champs de Markov gaussiens comme des ensembles d'approximation de la distribution de $X$, nous introduisons une nouvelle méthode de sélection de modèle pour des champs spatiaux. Pour tout voisinage $m$ dans une collection $\mathcal{M}$ donnée, cette procédure estime la covariance de $X$ par un champ de Markov de voisinage $m$. Puis, elle sélectionne un voisinage $\widehat{m}$ grâce à une technique de pénalisation. L'estimateur ainsi défini satisfait une inégalité oracle non-asymptotique. Si $X$ est un champ de Markov gaussien, la procédure est minimax adaptative à  la taille de son voisinage. Plus généralement, nous prouvons que la procédure s'adapte à  la vitesse d'approximation de la distribution de $X$ par des champs de Markov gaussiens de voisinage croissant.
}

\RRabstract{We study the nonparametric covariance estimation of a stationary Gaussian field $X$ observed on a regular lattice. In the time series setting, some procedures like AIC are proved to achieve optimal model selection among autoregressive models. However, there exists no such equivalent results of adaptivity in a spatial setting.  
By considering collections of Gaussian  Markov random fields (GMRF) as approximation sets for the distribution of $X$, we introduce a novel model selection procedure for spatial fields.
For all neighborhoods $m$ in a given collection $\mathcal{M}$, this procedure first amounts to computing a covariance  estimator of $X$ within the  GMRFs of neighborhood $m$. Then, it selects a neighborhood $\widehat{m}$ by applying a penalization strategy. The so-defined method satisfies a nonasymptotic oracle type inequality. If $X$ is a GMRF, the procedure is also minimax adaptive to the sparsity of its neighborhood. More generally, the procedure is adaptive to the rate of approximation of the true distribution by GMRFs with growing neighborhoods.
}

\RRmotcle{Champ gaussien, champ de Markov gaussien, sélection de modèle, pseudo-vraisemblance, inégalités oracles, vitesse minimax d'estimation.}
\RRkeyword{Gaussian field, Gaussian Markov random field, model selection, pseudolikelihood, oracle inequalities, Minimax rate of estimation.}

\RRprojets{Select}
\RRtheme{\THCog} % cas de 5 themes
 \RCSaclay % Saclay \^Ile de France
\begin{document}
\RRNo{6797} 
\makeRR   % cas d'un rapport de recherche

\section{Introduction}\label{introduction}

In this paper, we study the estimation of the distribution of a stationary Gaussian field $X=(X{\scriptstyle[i,j]})_{(i,j)\in\Lambda}$ indexed by the nodes of a square lattice $\Lambda$ of size $p\times p$. This problem is often encountered in spatial statistics or in image analysis.\\

% Parler des methodes d'estimation de la covariance

Various estimation methods have been proposed to handle this question. Most of them fall into two categories. On the one hand, one may consider direct covariance estimation.  A traditional approach amounts to computing an empirical variogram and then fitting a suitable parametric variogram model such as the exponential or Mat\'ern model (Cressie \cite{cressie} Ch.2). Some procedures also apply to non-regular lattices. However, a bad choice of the variogram model may lead to poor results. The issue of variogram model selection has not been completely solved yet, although some procedures based on cross-validation have been proposed. See \cite{cressie} Sect.2.6.4 for a discussion.
Most of the nonparametric~ (Hall \emph{et al.}~\cite{hall94}) and semiparametric (Im\emph{ et al.}~\cite{stein07}) methods are based on the spectral representation of the field. To our knowledge, these procedures have not yet been shown to achieve adaptiveness, i.e. their rate of convergence does not adapt to the \emph{complexity} of the correlation functions. \\

% Parler des methodes d'estimation type GMRF

An alternative approach to the problem amounts to considering the conditional distribution at one node given the remaining nodes. This point of view is closely connected to the notion of \emph{Gaussian Markov Random field} (GMRF). Let $\mathcal{G}$ be a graph whose vertex set is $\Lambda$. The field $X$ is GMRF with respect to $\mathcal{G}$ if it satisfies the following property: for any node $(i,j)\in\Lambda$, conditionally to the set of variables $X{\scriptstyle[k,l]}$ such that  $(k,l)$ is a neighbor of $(i,j)$ in $\mathcal{G}$,  $X{\scriptstyle[i,j]}$ is independent from all the remaining variables.  GMRFs are also  sometimes called Gaussian graphical models. A huge literature develops around this subject since Gaussian graphical models are  promising tools to analyze complex high-dimensional systems involved for instance in postgenomic data. In other applications, GMRFs are relevant because they allow to perform Markov chain Monte Carlo run fastly using Markov properties (e.g. \cite{rue02}). See Lauritzen~\cite{lauritzen} or  Edwards~\cite{Edwards} for introductions to Gaussian graphical models and Markov properties. In the sequel, we assume that the node $(0,0)$ belongs to $\Lambda$. Since we assume that the field $X$ is stationary, defining a graph $\mathcal{G}$ is equivalent to defining the neighborhood $m$ of the node $(0,0)$. Indeed, the neighborhood of any node $(i,j)\in\Lambda$ is the transposition of $m$ by $(i,j)$. In the sequel, we call $m$  \emph{the neighborhood} of a GMRF. If the neighborhood is empty, then the Markov property states that the components of $X$ are all independent. Alternatively, any zero-mean Gaussian stationary field is a GMRF with respect to the complete neighborhood (i.e. containing all the nodes except $(0,0)$).

Numerous papers have been devoted to parametric estimation for stationary GMRFs with a known neighborhood.  The authors have derived their asymptotic properties of such estimators (see \citep{besag75,besag77,Guyon}). If the field $X$ is assumed to be a GMRF with respect to a \emph{known} neighborhood in all these works, the issue of neighborhood selection has been less studied.  Besag and Kooperberg \cite{besag95}, Rue and Tjelmeland \cite{rue02}, Song \emph{et al.} \cite{fuentes2008}, and Cressie and Verzelen \cite{verzelen_cressie} have tackled  the problem of \emph{approximating} the distribution of a Gaussian field by a GMRF, but this requires the knowledge of the true distribution. Guyon and Yao have stated in \cite{guyon2000}  necessary conditions and sufficient conditions for a model selection procedure to choose asymptotically the true neighborhood of a GMRF with probability one.  \\

%%%%%%%%%%%%%%%%%%%%%%%%%%%%%%%%%%%%%%%%%%%%%%%%%%%%%%%%%%%%%%%%%%%%%%%%%%%%%%%%%
%%%%%%%%%%%%%%%%%%%%%%%%%%%%%%%%%%%%%%%%%%%%%%%%%%%%%%%%%%%%%%%%%%%%%%%%%%%%%%%%%
%%%%%%%%%%% INTRODUIRE LA PROBLEMATIQUE. %%%%%%%%%%%%%%%%%%%%%%%%%%%%%%%%%%%%%%%%
%%%%%%%%%%%%%%%%%%%%%%%%%%%%%%%%%%%%%%%%%%%%%%%%%%%%%%%%%%%%%%%%%%%%%%%%%%%%%%%%%
%%%%%%%%%%%%%%%%%%%%%%%%%%%%%%%%%%%%%%%%%%%%%%%%%%%%%%%%%%%%%%%%%%%%%%%%%%%%%%%%%

In this paper, we study a nonparametric estimation procedure based on neighborhood selection. 
In short, we select a \emph{suitable} neighborhood and estimate the distribution of $X$ in the space of stationary GMRFs with respect to this neighborhood. The objective is not to estimate the ``true'' neighborhood. We rather want to select a neighborhood that allows to estimate \emph{well} the distribution of $X$ (i.e. to minimize a risk). 
In fact, we do not even assume that the true correlation of $X$ corresponds to a GMRF. This estimation procedure is relevant for two main reasons:
\begin{itemize}		
 \item To our knowledge, it is the first nonparametric estimator in a spatial setting which achieves adaptive rates of convergence.
\item In most of the statistical applications where GMRFs are involved, the neighborhood is a priori unknown. Our procedure allows to select a ``good'' neighborhood.
\end{itemize}

%%%%%%%%%%%%%%%%%%%%%%%%%%%%%%%%%%%%%%%%%%%%%%%%%%%%%%%%%%%%%

%% Parler des liens avec les series temporelles et dans des cadres differents...

Our problem on a two-dimensional field has a natural one-dimensional counterpart in time series analysis. It is indeed known that an auto-regressive process (AR) of order $p$  is also a GMRF with $2p$ nearest neighbors and reciprocally (see \cite{guyon95} Sect. 1.3). In this one-dimensional setting, our issue reformulates as follows: how can we select the order of an AR to 
estimate well the distribution of a time series? It is known that  order selection by minimization of criteria like AICC, AIC or FPE satisfy asymptotically oracle inequalities (Shibata \cite{shibata80} and Hurvich and Tsai \cite{hurvich89}). We refer to Brockwell and Davis \cite{brockwell} and McQuarrie and Tsai \cite{MacQuarrie} for detailed discussions. However, one cannot readily extend these results to a spatial setting
because of computational and theoretical difficulties. 

In the rest of this introduction, we further describe the framework and we summarize the main results of the paper.

\subsection{Conditional regression}\label{section_approche}

Let us now make precise the notations and present the ideas underlying our approach. 
In the sequel,  $\Lambda$ stands for the  toroidal lattice of size $p\times p$. We consider the random field  $X=(X{\scriptstyle [i,j]})_{1\leq i,j\leq p}$  indexed by the nodes of $\Lambda$. Besides, $X^v$ refers  to the vectorialized version of $X$ with the convention $X{\scriptstyle[i,j]}=X^v{\scriptstyle[(i-1)\times p+ j]}$ for any $1\leq i,j\leq p$. Using this new notation amounts to ``forgetting'' the spatial structure of $X$ and allows to get into a more classical statistical framework. For the sake of simplicity, the components of $X$ are defined modulo $p$ in the remainder of the paper.

Throughout this paper, we assume the field $X$ is centered. In practice, the statistician has to first subtract some parametric form of the mean value. Hence, the vector $X^v$ follows a zero-mean Gaussian distribution $\mathcal{N}(0,\Sigma)$,  where the $p^2\times p^2$ matrix $\Sigma$ is non singular but unknown. Besides, we suppose that the field $X$ is stationary on the torus $\Lambda$. More precisely, for any $r>0$,
any $(i,j)\in \{1,\ldots,p\}^2$, and any $(k_1,l_1),\ldots,(k_r,l_r)\in
\{1,\ldots,p\}^{2r}$, it holds that
$$\left(X{\scriptstyle[k_1,l_1]},\ldots,X{\scriptstyle[k_r,l_r]}\right)\sim \left(X{\scriptstyle[k_1+i,l_1+j]},\ldots,X{\scriptstyle[k_r+i,l_r+j]}\right) \ . $$

We observe $n\geq 1$ i.i.d. replications of the vector $X^v$. In the sequel, ${\bf X^v}$ denotes the $p^2\times n$ matrix of the $n$ observations of $X^v$. For any $1\leq i\leq n$, the $p\times p$ matrix ${\bf X}_i$ stands for the $i$-th observation of the field $X$. All these notations are recalled in Table \ref{tableau_notation} in Section \ref{section_notations}.
In practice, the number of observations $n$ often equals one. Our goal is to estimate the matrix $\Sigma$.\\

We sometimes assume that the field $X$ is isotropic. Let $G$ be the group of vector isometries of the unit square. For any node $(i,j)\in\Lambda$ and any isometry $g\in G$, $g.(i,j)$ stands for the image of $(i,j)$ in $\Lambda$ under the action of $g$. We say that $X$ is isotropic on $\Lambda$ if for any $r>0$,
$g\in G$, and $(k_1,l_1),\ldots,(k_r,l_r)\in
\{1,\ldots,p\}^{2r}$, 
$$\left(X{\scriptstyle[k_1,l_1]},\ldots,X{\scriptstyle[k_r,l_r]}\right)\sim \left(X{\scriptstyle[g.(k_1,l_1)]},\ldots,X{\scriptstyle[g.(k_r,l_r)]}\right) \ . $$

As mentioned earlier, we aim at estimating the distribution of the field $X$ through a conditional distribution approach. By standard Gaussian derivations (see for instance \cite{lauritzen} App.C), there exists a unique $p\times p$ matrix $\theta$ such that $\theta{\scriptstyle[0,0]}=0$ and
\begin{eqnarray}\label{regression_conditionnelle}
X{\scriptstyle[0,0]}= \sum_{(i,j)\in \Lambda\setminus\{(0,0)\}}\theta{\scriptstyle[i,j]}X{\scriptstyle[i,j]} + \epsilon{\scriptstyle[0,0]}\ ,
\end{eqnarray}
where the random variable $\epsilon{\scriptstyle[0,0]}$ follows a zero-mean normal distribution and is independent from the covariates $(X{\scriptstyle[i,j]})_{(i,j)\in \Lambda\setminus\{(0,0)\}}$. Equation (\ref{regression_conditionnelle}) describes the conditional distribution of $X{\scriptstyle[0,0]}$ given the remaining variables. Since the field $X$ is stationary, the matrix $\theta$ also satisfies $\theta{\scriptstyle[i,j]} = \theta{\scriptstyle[-i,-j]}$ for any $(i,j)\in \Lambda$.
Let us note $\sigma^2$ the conditional variance of $X{\scriptstyle[0,0]}$ and $I_{p^2}$ the identity matrix of size $p^2$. The matrix $\theta$ is closely related to the covariance matrix $\Sigma$ of $X^v$ through the following property: 
\begin{eqnarray}\label{lien_sigma_theta}
 \Sigma = \sigma^2 \left(I_{p^2} - C(\theta)\right)^{-1} \ ,
\end{eqnarray}
where the $p^2\times p^2$ matrix $C(\theta)$ is defined as $C(\theta){\scriptstyle[i_1(p-1)+j_1,i_2(p-1)+j_2]} := \theta{\scriptstyle[i_2-i_1,j_2-j_1]}$ for any $1\leq i_1,i_2,j_1,j_2\leq p$. The matrix $(I_{p^2} - C(\theta))$ is called the partial correlation matrix of the field $X$.
The so-defined matrix $C(\theta)$ is symmetric block circulant with $p\times p$ blocks as stated below. We refer to \cite{rue} Sect.2.6 or the book of Gray \cite{gray} for definitions and main properties on circulant and block circulant matrices.
 
\begin{lemma}\label{block_circulant}
Let $\theta$ be a square matrix of size $p$ such that 
\begin{eqnarray}\label{symmetric_matrice}
\text{for any }1\leq i,j\leq p,\, \theta{\scriptstyle[i,j]} = \theta{\scriptstyle[-i,-j]},
\end{eqnarray}
then the matrix $C(\theta)$ is symmetric block circulant with $p\times p$
blocks. Conversely, if $B$ is a $p^2\times p^2$ symmetric block circulant matrix
with $p\times p$ blocks, then there exists a square matrix $\theta$ of size
$p$ satisfying (\ref{symmetric_matrice}) and such that $B=C(\theta)$.
\end{lemma}
A proof is given in the technical appendix~\cite{technical}. In conclusion, estimating the matrix $\Sigma/\sigma^2$ amounts to estimating the matrix $C(\theta)$, which is also equivalent to estimating the $p\times p$ matrix $\theta$. This is why, we shall focus on the estimation of the matrix $\theta$. \\

Let us precise the set of possible values for $\theta$. In the sequel,  $\Theta$ denote the vector space of the $p\times p$ matrices that satisfy $\theta{\scriptstyle[0,0]}=0$ and 
$\theta{\scriptstyle[i,j]}=\theta{\scriptstyle[-i,-j]}$, for any $(i,j)\in\Lambda$. A matrix $\theta\in\Theta$ corresponds to the distribution of a stationary Gaussian field if and only if the $p^2\times p^2$ matrix $(I_{p^2}-C(\theta))$ is positive definite. This is why we define the convex subset $\Theta^+$ of $\Theta$ by
\begin{eqnarray}
 \Theta^+:=\left\{\theta\in\Theta\hspace{0.2cm}\text{s.t.}\,  \left(I_{p^2}-C(\theta)\right)\text{ is positive definite}\right\}\label{definition_theta+}\ .
\end{eqnarray}
The set of covariance matrices of stationary Gaussian fields on $\Lambda$ with unit conditional variance is therefore in one to one correspondence with the set $\Theta^+$.  Let us define the corresponding set $\Theta^{\text{iso}}$ and $\Theta^{+,\text{iso}}$ for isotropic Gaussian fields.
\begin{eqnarray}\label{def_Thetaiso}
 \Theta^{\text{iso}} :=\left\{\theta\in\Theta\ ,\theta{\scriptstyle[i,j]}=\theta{\scriptstyle[g.(i,j)]}\ ,\,  \forall (i,j)\in\Lambda,\ \forall g\in G \right\}\text{ and } \Theta^{+,\text{iso}} := \Theta^+\cap\Theta^{\text{iso}}\ .
\end{eqnarray}

\subsection{Model selection}

We have the issue of covariance estimation as an estimation problem for conditional regressions (Equation (\ref{regression_conditionnelle})). However, the set $\Theta^+$ of admissible parameters for the estimation is huge. The dimension of $\Theta$ is indeed of the same order as $p^2$ whereas we only observe $p^2$ non-independent data if $n$ equals one. In order to avoid the curse of dimensionality, it is natural to assume that the target $\theta$ is approximately \emph{sparse}. \\

It is indeed likely that the coefficients $\theta{\scriptstyle[i,j]}$ are \emph{close} to zero for the nodes $(i,j)$ which are \emph{far} from the origin $(0,0)$. By Equation (\ref{regression_conditionnelle}), this means that $X{\scriptstyle[0,0]}$ is \emph{well} predicted by the covariates $X{\scriptstyle[i,j]}$ whose corresponding nodes $(i,j)$ are close to the origin. In other terms, the true covariance is presumably well approximated by a GMRF with a \emph{reasonable} neighborhood. The main difficulty is that we do not know a priori what ``reasonable'' means. We want to adapt to the \emph{sparsity} of the matrix $\theta$. \\

% Je reecris le paragraphe précédent en fusionnant les deux chapitres et en écrivant un seul...
In the sequel, $m$ refers to a subset of $\Lambda\setminus\{0,0\}$. We call it a model. By Equation (\ref{regression_conditionnelle}),
the property ``$X$ is a GMRF with respect to the neighborhood $m$'' is  equivalent to ``the support of $\theta$ is included in $m$''.
We are given a nested collection $\mathcal{M}$ of models. For any of these models $m\in\mathcal{M}$, we compute $\widehat{\theta}_{m,\rho_1}$ the Conditional least squares estimator (CLS) of $\theta$ for the model $m$ by maximizing the pseudolikelihood over a subset of matrices $\theta$ whose support is included in $m$. These estimators as well as their dependency on the quantity $\rho_1$ are defined in Section \ref{section_procedure}. \\

The model $m$ that minimizes the risk of $\widehat{\theta}_{m,\rho_1}$ over the collection $\mathcal{M}$ is called an oracle and is noted $m^*$. In practice, this model is unknown and we have to estimate it. The art of model selection is to pick a model $m\in\mathcal{M}$ that is large enough to enable a good approximation of $\theta$ but is small enough so that the variance of $\widehat{\theta}_{m,\rho_1}$ is small. Let us reformulate the approach in terms of GMRFs: given a collection $\mathcal{M}$ of neighborhoods, we compute an estimator of $\theta$ in the set of GMRFs with neighborhood $m$,  for any $m\in\mathcal{M}$. Our purpose is to select a suitable neighborhood $\widehat{m}$ so that the estimator $\widehat{\theta}_{\widehat{m}}$ has a risk as small as possible.

A classical method to estimate a \emph{good model} $\widehat{m}$ is achieved through \emph{penalization} with respect to the size of the models. In the following expression, $\gamma_{n,p}(.)$ stands for the CLS empirical contrast that we shall define in Section \ref{section_procedure}.  We select a model $\widehat{m}$ by minimizing the criterion
\begin{eqnarray}\label{definition_mhat}
\widehat{m} = \arg\min_{m\in\mathcal{M}} \left[\gamma_{n,p}(\widehat{\theta}_{m,\rho_1})+\pen(m)\right]\ .
\end{eqnarray}
where $\pen(.)$ denotes a positive function defined on $\mathcal{M}$.
In this paper, we prove that under a suitable choice of the penalty function $\pen(.)$, the risk of the estimator $\widehat{\theta}_{\widehat{m}}$ is as small as possible.

\subsection{Risk bounds and adaptation}

We shall assess our procedure using two different loss functions. First, we introduce the loss function $l(.,.)$ that measures how well we estimate the conditional distribution (\ref{regression_conditionnelle}) of the field. For any $\theta_1,\theta_2\in\Theta$, the  distance $l(\theta_1,\theta_2)$ is defined by
\begin{eqnarray}\label{def_lmatriciel}\label{definition_l}
l\left(\theta_1,\theta_2\right) &:= &\frac{1}{p^2}tr\left[\left(C(\theta_1)-C(\theta_2)\right)\Sigma
  \left(C(\theta_1)-C(\theta_2)\right)\right ]\ .
\end{eqnarray}
Let us reformulate $l(\theta_1,\theta_2)$ in terms of conditional expectation
\begin{eqnarray*}
l\left(\theta_1,\theta_2\right)=  \mathbb{E}_{\theta}\left\{\left[\mathbb{E}_{\theta_1}\left(X{\scriptstyle[0,0]}|X_{\Lambda\setminus\{0,0\}}\right)-\mathbb{E}_{\theta_2}\left(X{\scriptstyle[0,0]}|X_{\Lambda\setminus\{0,0\}}\right)\right]^2\right\}\ ,
\end{eqnarray*}
where $\mathbb{E}_{\theta}(.)$ stands for the expectation with respect to the distribution of $X^v$,  $\mathcal{N}(0,\sigma^2(I_{p^2}-C(\theta))^{-1})$. Hence, $l(\widehat{\theta},\theta)$ corresponds the mean squared prediction loss which is often used in the random design regression framework, in time series analysis \cite{hurvich89}, or in spatial statistics \cite{fuentes2008}.
Moreover, the loss function $l(\widehat{\theta},\theta)$ is also connected to the notion of kriging error. The kriging predictor (Stein~\cite{stein}) of $X{\scriptstyle[0,0]}$ is defined as the best linear combination of the covariates $(X{\scriptstyle [k,l]})_{(k,l)\in\Lambda\setminus\{(0,0\}}$ for predicting the value $X{\scriptstyle [0,0]}$. 
By Equation (\ref{regression_conditionnelle}), this predictor is exactly $\sum_{(k,l)\in\Lambda\setminus\{(0,0\}}\theta{\scriptstyle [k,l]}X{\scriptstyle [k,l]}$ and the mean squared prediction error is $\sigma^2$. If we do not know $\theta$ but we are given  an estimator $\widehat{\theta}$, then the corresponding kriging predictor $\sum_{(k,l)\in\Lambda\setminus\{(0,0\}}\widehat{\theta}{\scriptstyle [k,l]}X{\scriptstyle [k,l]}$ has a mean squared prediction error equal to $\sigma^2+l(\widehat{\theta},\theta)$. Kriging is a key concept in spatial statistics and it is therefore interesting to consider a loss function that measures the kriging performances when one estimates $\theta$.\\

 We shall also assess our results using the Frobenius distance noted $\|.\|_F$ and defined by $\|A\|_F^2:=\sum_{1\leq i,j\leq p}A{\scriptstyle [i,j]}^2$. Observe that the Frobenius distance $\|\theta_1-\theta_2\|_F^2$  also equals the Frobenius distance between the partial correlation matrices $(I_{p^2}-C(\theta_1))$ and $(I_{p^2}-C(\theta_2))$ (up to a factor $p^2$)
\begin{eqnarray}\label{perte_frobenius}
 \|\theta_1-\theta_2\|_F^2=\frac{1}{p^2}\|\big(I_{p^2}-C(\theta_1)\big)-\big(I_{p^2}-C(\theta_2)\big)\|_F^2\ ,
\end{eqnarray}

Our aim is then to define a suitable penalty function $\pen(.)$ in (\ref{definition_mhat}) so that the estimator $\widehat{\theta}_{\widehat{m},\rho_1}$ 
performs almost as well as the oracle estimator $\widehat{\theta}_{m^*,\rho_1}$. For any model $m\in\mathcal{M}$, we define $\theta_{m,\rho_1}$ as the matrix which minimizes the loss $l(\theta',\theta)$ over the sets of matrices $\theta'$ corresponding to  model $m$. The loss $l(\theta_{m,\rho_1},\theta)$ is called the \emph{bias}. Our main result is stated in Section \ref{section_main_result}. We provide a condition on the penalty function $\pen(.)$, so that the selected estimator satisfies a risk bound of the form
\begin{eqnarray}\label{inegalite_oracle_introduction}
 \mathbb{E}_{\theta}\left[l\left(\widehat{\theta}_{\widehat{m},\rho_1},\theta\right)\right]\leq L\inf_{m\in\mathcal{M}}\left[l(\theta_{m,\rho_1},\theta)+\varphi_{\text{max}}(\Sigma)\frac{\text{Card}(m)}{np^2}\right] \ ,
\end{eqnarray}
where $\varphi_{\text{max}}(\Sigma)$ is the largest eigenvalue of $\Sigma$ and $\text{Card}(.)$ stands for the cardinality.
Contrary to most results in a spatial setting, this upper bound on the risk is nonasymptotic and holds in a general setting. The term $\varphi_{\text{max}}(\Sigma)\text{Card}(m)/(np^2)$ grows linearly with the size of $m$ and goes to $0$ with $n$ and $p$. In Section \ref{section_risk}, we prove that the variance term of a model $m$ is of the same order as $\varphi_{\text{max}}(\Sigma)\text{Card}(m)/(np^2)$. Hence, the bound (\ref{inegalite_oracle_introduction}) tells us that the risk of $\widehat{\theta}_{\widehat{m},\rho_1}$ is smaller than a quantity which is the same order as the risk $\mathbb{E}_{\theta}[l(\widehat{\theta}_{m^*,\rho_1},\theta)]$ of the oracle $m^*$. We say that the selected estimator achieves an \emph{oracle-type inequality}. \\

In Section \ref{section_risk}, we bound the asymptotic expectations $\mathbb{E}[l(\widehat{\theta}_{m,\rho_1},\theta)]$ and connect them to the variance terms in Bound (\ref{inegalite_oracle_introduction}). As a consequence, we prove that under mild assumptions on the target $\theta$, the upper bound (\ref{inegalite_oracle_introduction}) is optimal from the asymptotic point of view (up to a multiplicative numerical constant). We discuss the assumptions in Section \ref{section_discussion_hypotheses}.
In Section \ref{section_minimax}, we compute nonasymptotic minimax lower bounds with respect to the loss functions $l(.,.)$ and $\|.\|_F^2$. We then derive that under mild assumptions, our estimator $\widehat{\theta}_{\widehat{m},\rho_1}$ is minimax adaptive to the sparsity of $\theta$ and minimax adaptive to the decay of $\theta$.\\

To our knowledge, these are the first oracle-type inequalities in a spatial setting.
The computation of the minimax rates of convergence is also new. Moreover, most of our results are nonasymptotic. Although we have considered a square on the two-dimensional lattice, our method straightforwardly extends to any $d$-dimensional toroidal rectangle with $d\geq 1$. In the one-dimensional setting, we retrieve a oracle-type inequality that is close to the work of  Shibata \cite{shibata80}. Yet, he has stated an asymptotic oracle inequality for the estimation of autoregressive processes. In contrast, our result applies on a torus and is only optimal up to constants but it is nonasympotic and most of all applies for higher dimensional lattices.
In Section \ref{section_discussion}, we further discuss the advantages and the weak points of our method. Moreover, we mention the extensions and the simulations made in a subsequent paper \cite{verzelen_gmrf_appli}. All the proofs are postponed to Section \ref{section_Preuves} and to the appendix~\cite{technical}.

\subsection{Some notations}\label{section_notations}

Throughout this paper, $L,L_1, L_2, \ldots$ denote constants that may vary from line to line. The notation $L(.)$ specifies the dependency on some quantities. For any matrix $A$,  $\varphi_{\text{max}}(A)$ and  $\varphi_{\text{min}}(A)$ respectively refer the largest eigenvalue and the smallest eigenvalues of $A$. We recall that $\|A\|_F$ is the Frobenius norm of $A$.
For any matrix $\theta$ of size $p$, $\|\theta\|_1$ stands for the sum of of the absolute values of the components of $\theta$, we call it its $l_1$ norm. In the sequel, $0_p$ is the square matrix of size $p$ whose indices are $0$. Given $\rho>0$, the ball $\mathcal{B}_1(0_p;\rho)$ is defined as the set of square matrices of size $p$ whose $l_1$ norm is smaller than $\rho$. Finally, Table \ref{tableau_notation} gathers the notations involving $X$.

\begin{Table}[h]
\caption{Notations for the random field and the data.\label{tableau_notation}}
\begin{center}
\begin{tabular}{|c|l|l|l}
 \hline $X$ & Matrix of size $p\times p$ & Random field\\ \hline
$X^v$ & Vector of length $p^2$ & Vectorialized version of $X$ \\ \hline
${\bf X^v}$ & Matrix of size $p^2\times n$ & Observations of $X^v$\\ \hline 
${\bf X}_i$ & Matrix of size $p\times p$ & $i$-th observation of the field $X$\\ \hline
\end{tabular}
\end{center}
\end{Table}

\section{Model selection procedure }\label{section_procedure}

In this section, we formally define our model selection procedure.

\subsection{Collection of models}

For any node $(i,j)$ belonging to the lattice $\Lambda$, let us define the toroidal norm by
\begin{eqnarray*}
|(i,j)|^2_t :=\left[i\wedge \left(p-i\right)\right]^2 +\left[j\wedge \left(p-j\right)\right]^2
\end{eqnarray*}

We aim at selecting a ``good'' neighborhood for the GMRF. Since $X$ corresponds to some ``spatial'' process, it is natural to assume that nodes that are close to $(0,0)$ are more likely to be significant. This is why we restrict ourselves in the sequel to the collection $\mathcal{M}_1$ of neighborhoods.

\begin{defi}\label{definition_modele}
A subset $m\subset \Lambda\setminus\{(0,0)\}$ belongs to $\mathcal{M}_{1}$ if  
 there exists a number $r_m>1$ such that 	
\begin{eqnarray}
m = \left\{ (i,j)\in \Lambda\setminus\{(0,0)\}\hspace{0.3cm} \text{s.t.}\hspace{0.3cm}  |(i,j)|_t\leq r_m\right\}\ . \label{definition_rm}
\end{eqnarray}
\end{defi}

\begin{figure}
\centerline{\epsfig{file=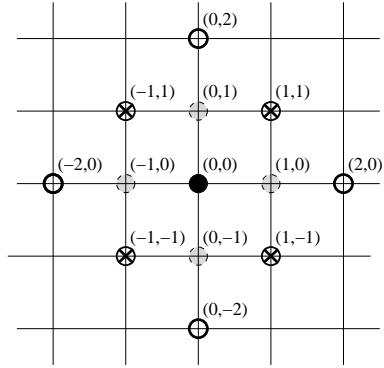,
    width=5cm}}
\caption{\textit{Examples of models. The four gray nodes refer to $m_1$. The model $m_2$ also contains the nodes with a cross whereas $m_3$ contains all the nodes except $(0,0)$. }\label{exemple-modeles}}
\end{figure}

The collection $\mathcal{M}_{1}$ is totally ordered with respect to the
inclusion and we therefore order our models $m_0\subset m_1\subset \ldots \subset m_{i}\ldots  $. For instance, $m_0$ corresponds to the empty neighborhood whereas $m_1$ stands for the neighborhood of size $4$. See Figure \ref{exemple-modeles} for other examples.\\

For any model $m\in\mathcal{M}_1$, we define the vector space $\Theta_m$ as the subset of the elements of $\Theta$ whose support is included in $m$. We recall that $\Theta$ is defined in Section \ref{section_approche}. Similarly $\Theta^{\text{iso}}_m$ is the subset of $\Theta^{\text{iso}}$  whose support is included in $m$. The dimensions of $\Theta_m$  and $\Theta^{\text{iso}}_m$ are respectively noted $d_m$ and $d_m^{\text{iso}}$. 
Since we aim at estimating the positive matrix $(I_{p^2}-C(\theta))$, we shall consider the convex subsets
 of $\Theta_m^+$ and $\Theta_m^{+,\text{iso}}$ that correspond to non-negative precision matrices.
\begin{eqnarray}\label{definition_thetam+}
 \Theta_m^+ :=\Theta_m\cap \Theta^+ \hspace{1.5cm}\text{ and }\hspace{1.5cm} \Theta_m^{+,\text{iso}} :=\Theta^{\text{iso}}_m\cap \Theta^{+,\text{iso}}\ .
\end{eqnarray}
For instance, the set $\Theta_{m_1}^+$ is in one to one correspondence with the sets of GMRFs whose neighborhood is made of the four nearest neighbors. Similarly, $\Theta_{m_1}^+$ is in one to one correspondence with the GMRFs with eight nearest neighbors.
In our estimation procedure, we shall restrict ourselves to precision matrices whose largest eigenvalue is upper bounded by a constant. This is why we define the subsets $\Theta_{m_2,\rho_1}^{+}$ and $\Theta_{m,\rho_1}^{+,\text{iso}}$ for any  $\rho_1\geq 2$.
\begin{eqnarray}
\Theta_{m,\rho_1}^{+} & := &\left\{\theta\in  \Theta_m^{+}\ , \varphi_{\text{max}}\left(I_{p^2}-C(\theta)\right)< \rho_1\right\} \label{defintion_thetam+gamma}\\
\Theta_{m,\rho_1}^{+,\text{iso}}  & := &\left\{\theta\in  \Theta_m^{+,\text{iso}}\ , \varphi_{\text{max}}\left(I_{p^2}-C(\theta)\right)< \rho_1\right\} \ . \label{defintion_thetam+gamma_iso}
\end{eqnarray}\hspace{1.5cm}

Finally, we need a generating family of the spaces $\Theta_m$ and $\Theta_m^{\text{iso}}$. 
For any node $(i,j)\in  \Lambda\setminus \{(0,0)\}$, let us define the $p\times p$ matrix $\Psi_{i,j}$ as 
\begin{eqnarray}
\Psi_{i,j}{\scriptstyle[k,l]}:= \left\{\begin{array}{ccc} 1 & \text{if}&
\text{$(k,l)=(i,j)$ or  $(k,l)=-(i,j)$}\\ 0 & \text{otherwise\ .}&\end{array} \right. \label{definition_psi}
\end{eqnarray}
Hence,  $\Theta_m$ is generated by the matrices $\Psi_{i,j}$ for which $(i,j)$ belongs to $m$. Similarly, for any $(i,j)\in \Lambda\setminus \{(0,0)\}$, let us define the matrix $\Psi_{i,j}^{\text{iso}}$ by
\begin{eqnarray}\label{definition_psiiso}
\Psi_{i,j}^{\text{iso}}{\scriptstyle[k,l]}:= \left\{\begin{array}{ccc} 1 & \text{if} &
  \text{$\exists g\in G$, $(k,l)=g.(i,j)$}\\ 0 & \text{otherwise\ .}&\end{array} \right.
\end{eqnarray}

\subsection{Estimation by Conditional Least Squares (CLS)}\label{section_cls}

Let us turn to  the conditional least squares estimator. For any $\theta'\in \Theta^+$, the  criterion $\gamma_{n,p}(\theta')$ is defined by
\begin{eqnarray}\label{definition_gamman}
\gamma_{n,p}(\theta') & := & \frac{1}{np^2}\sum_{i=1}^n\sum_{1\leq j_1,j_2\leq p} \bigg({\bf X}_i{\scriptstyle[j_1,j_2]} -
\sum_{(l_1,l_2)\in \Lambda\setminus\{(0,0)\}} \theta'{\scriptstyle[l_1,l_2]}{\bf X}_i{\scriptstyle[j_1+l_1,j_2+l_2]}\bigg)^2 \ .
\end{eqnarray}
In a nutshell, $\gamma_{n,p}(\theta')$ is a least squares criterion that allows to perform the simultaneous linear regression of all ${\bf X}_i{\scriptstyle[j_1,j_2]}$ with respect to the covariates $({\bf X}_i{\scriptstyle[l_1,l_2]})_{(l_1,l_2)\neq(j_1,j_2)}$. The advantage of this criterion is that it does not require the computation of a determinant of a huge matrix as for the likelihood.
We shall often use an alternative expression of $\gamma_{n,p}(\theta')$ in terms of the factor $C(\theta')$ 
and the empirical covariance matrix $\overline{{\bf X^v X^{v*}}}$:
\begin{eqnarray}\label{definition_gamman_alternative}
\gamma_{n,p}(\theta')& = & \frac{1}{p^2}tr\left[(I_{p^2}-C(\theta'))\overline{{\bf X^v X^{v*}}}(I_{p^2}-C(\theta'))\right] \ . 
\end{eqnarray}
One proves the equivalence between these two expressions by coming back to the definition of $C(\theta')$.
Let $\rho_1>2$ be fixed. For any model $m\in\mathcal{M}$, we compute the CLS estimators $\widehat{\theta}_{m,\rho_1}$ and $\widehat{\theta}^{\text{iso}}_{m,\rho_1}$ by minimizing the criterion $\gamma_{n,p}(.)$ as follows
\begin{eqnarray}\label{definition_estimateur_cls}
 \widehat{\theta}_{m,\rho_1}:=\arg\min_{\theta'\in\overline{\Theta_{m,\rho_1}^+}} \gamma_{n,p}(\theta')\hspace{1cm}\text{and}\hspace{1cm}  \widehat{\theta}^\text{iso}_{m,\rho_1}:=\arg\min_{\theta'\in\overline{\Theta_{m,\rho_1}^{+,\text{iso}}}} \gamma_{n,p}(\theta')\ ,
\end{eqnarray}
where $\overline{A}$ stands for the closure of the set $A$. The existence and the uniqueness of $\widehat{\theta}_{m,\rho_1}$ and $\widehat{\theta}^\text{iso}_{m,\rho_1}$ are ensured  by the following lemma.
\begin{lemma}\label{gamma_convexe}
For any $\theta\in \Theta^+$,
$\gamma_{n,p}(.)$ is almost surely strictly convex on $\overline{\Theta^+}$.
\end{lemma}
The proof is postponed to the appendix~\cite{technical}. We discuss the dependency of $\widehat{\theta}_{m,\rho_1}$ on  the parameter $\rho_1$ in Section \ref{section_discussion_hypotheses}.
For stationary Gaussian fields, minimizing the CLS criterion $\gamma_{n,p}(.)$ over a  set $\Theta_{m,\rho_1}^+$
 is equivalent to minimizing the product of the conditional likelihoods $(X{\scriptstyle[i,j]}|X_{-\{i,j\}})$, called \emph{Conditional Pseudo-Likelihood} (CPL): 
\begin{eqnarray*}
 p\mathcal{L}_n(\theta',{\bf X^v}) := \prod_{\begin{array}{c}
1\leq i\leq n,\\ (j_1,j_2)\in \Lambda\end{array}}\mathcal{L}_{n,\theta'}\left({\bf X}_i{\scriptstyle[j_1,j_2]}|({\bf X}_i)_{-\{j_1,j_2\}}\right)= \left(\sqrt{2\pi}\sigma\right)^{-np^2}\exp\left(- \frac{1}{2}\frac{np^2\gamma_{n,p}(\theta')}{\sigma^2}\right) \ , 
\end{eqnarray*}
where we recall that $\sigma^2$ refers to the conditional variance of any $X{\scriptstyle[i,j]}$. In fact, CLS estimators were first introduced by Besag \cite{besag75-2} who call them pseudolikelihood estimators since they minimize the CPL. \\

Let us define the  function $\gamma(.)$ as an infinite sampled version of 
the CLS criterion $\gamma_{n,p}(.)$:
\begin{eqnarray}\label{definition_gamma}
 \gamma(\theta') := \mathbb{E}_{\theta}\left[\gamma_{n,p}(\theta')\right]=\mathbb{E}_{\theta}\bigg[\bigg(X{\scriptstyle[0,0]} - \sum_{(i,j)\neq (0,0)}\theta'{\scriptstyle[i,j]}X{\scriptstyle[i,j]}\bigg)^2\bigg]\ ,
\end{eqnarray}
for any $\theta',\theta\in \Theta^{+}$. The function $\gamma(\theta')$ measures the prediction error of $X{\scriptstyle[0,0]}$ if one uses $\sum_{(i,j)\neq (0,0)}\theta'{\scriptstyle[i,j]}X{\scriptstyle[i,j]}$ as a predictor. Moreover, it is a special case of the CMLS criterion introduced by Cressie and Verzelen in (Eq.10) of \cite{verzelen_cressie}  to approximate a Gaussian field by a GMRF. Hence, one may interpret the CLS criterion as a finite sampled version of their approximation method.
Observe that the function $\gamma(.)$ is minimized over $\Theta^+$  at
the point $\theta$ and that $\gamma(\theta) = \var_{\theta}(X{\scriptstyle[0,0]}\left|X_{-\{0,0\}} \right.)=\sigma^2$. Moreover,  the difference $\gamma(\theta')-\gamma(\theta)$ equals the loss $l(\theta',\theta)$ defined by (\ref{definition_l}).\\

For any model $m\in\mathcal{M}$, we introduce the projections $\theta_{m,\rho_1}$ and $\theta^{\text{iso}}_{m,\rho_1}$ as the best approximation of $\theta$ in $\overline{\Theta_{m,\rho_1}^+}$ and $\overline{\Theta_{m,\rho_1}^{+,\text{iso}}}$.
\begin{eqnarray}
\theta_{m,\rho_1}:=\arg\min_{\theta'\in\overline{\Theta_{m,\rho_1}^+}} l(\theta',\theta)\hspace{1cm}\text{and}\hspace{1cm}  \theta^\text{iso}_{m,\rho_1}:=\arg\min_{\theta'\in\overline{\Theta_{m,\rho_1}^{+,\text{iso}}}} l(\theta',\theta)\ .
\label{def_thetam}
\end{eqnarray}
Since $\gamma(.)$ is strictly convex on $\Theta^+$, the matrices $\theta_{m,\rho_1}$ and $\theta^{\text{iso}}_{m,\rho_1}$ are uniquely defined. By its definition (\ref{definition_l}), one may interpret $l(.,.)$ as an inner product on the space $\Theta$; therefore, the orthogonal projection of $\theta$ onto the convex closed set $\overline{\Theta_{m,\rho_1}^+}$ (resp. $\overline{\Theta_{m,\rho_1}^{+,\text{iso}}}$) with respect to $l(.,.)$  is $\theta_{m,\rho_1}$ (resp. $\theta^\text{iso}_{m,\rho_1}$). It then follows from a property of orthogonal projections that the loss of $\widehat{\theta}_{m,\rho_1}$ is upper bounded by
\begin{eqnarray}\label{decomposition_biais_variance_inegalite}
l(\widehat{\theta}_{m,\rho_1},\theta)\leq l(\theta_{m,\rho_1},\theta)+l(\widehat{\theta}_{m,\rho_1},\theta_{m,\rho_1})\ .
\end{eqnarray}
The first term $l(\theta_{m,\rho_1},\theta)$ accounts for the bias, whereas the second term $l(\widehat{\theta}_{m,\rho_1},\theta_{m,\rho_1})$ is a variance term. Observe that $\theta\in\Theta^{+}_m$ does not necessarily imply that the bias $l(\theta_{m,\rho_1},\theta)$ is null because in general $\overline{\Theta_{m}^+}\neq \overline{\Theta_{m,\rho_1}^+}$. This will be the case only  if $\theta$ satisfies the following hypothesis.
\begin{eqnarray}\label{def_h1}
(\mathbb{H}_1):\hspace{3.5cm} \varphi_{\text{max}}(I_	{p^2}-C(\theta))<\rho_1\ . \hspace{3.5cm}
\end{eqnarray}

Assumption $(\mathbb{H}_1)$ is necessary to ensure the existence of a model $m\in\mathcal{M}$ such that the bias is zero (i.e. $\theta_{m,\rho_1}=\theta$). By identity (\ref{lien_sigma_theta}), one observes that $(\mathbb{H}_1)$ is equivalent to a lower bound on the smallest eigenvalue of $\Sigma$, i.e. $\varphi_{\text{min}}(\Sigma)\leq \sigma^2/\rho_1$. We further discuss $(\mathbb{H}_1$) in Section \ref{section_discussion_hypotheses}.\\ 

For the sake of completeness, we recall the penalization criterion introduced in (\ref{definition_mhat}).
Given a subcollection of models $\mathcal{M}\subset\mathcal{M}_1$ and a positive function $\pen:\mathcal{M}\rightarrow \mathbb{R}^+$ that we call a penalty, we select a model as follows
\begin{eqnarray*}
 \widehat{m} := \arg\min_{m\in\mathcal{M}}\left[\gamma_{n,p}\left(\widehat{\theta}_{m,\rho_1}\right)\right]+\pen(m)\hspace{0.5cm}\text{and}\hspace{0.5cm}\widehat{m}^{\text{iso}} := \arg\min_{m\in\mathcal{M}}\left[\gamma_{n,p}\left(\widehat{\theta}^{\text{iso}}_{m,\rho_1}\right)\right]+\pen(m)\ .
\end{eqnarray*}
Observe that $\widehat{m}$ and $\widehat{m}^{\text{iso}}$ depend on $\rho_1$. For the sake clarity, we do not emphasize this dependency in the notation.
In the sequel, we write $\widetilde{\theta}_{\rho_1}$ and $\widetilde{\theta}_{\rho_1}^{\text{iso}}$ for $\widehat{\theta}_{\widehat{m},\rho_1}$ and $\widehat{\theta}_{\widehat{m}^{\text{iso}}}^{\text{iso},\rho_1}$.\\

\section{Main Result}\label{section_main_result}

We now provide a nonasymptotic upper bound for the risk of the estimators
 $\widetilde{\theta}_{\rho_1}$ and $\widetilde{\theta}_{\rho_1}^{\text{iso}}$. Let us recall that $\Sigma$ stands for the covariance matrix of $X^v$.

\begin{thrm}\label{mainthrm}
Let $K$ be a positive number larger than a universal constant $K_0$ and let $\mathcal{M}$ be a subcollection of $\mathcal{M}_1$. If for every model $m\in \mathcal{M}$, 
\begin{eqnarray} \label{condition_penalite}
\pen(m) \geq K\rho_1^2\varphi_{\text{max}}(\Sigma)\frac{d_m}{np^2}\ ,
\end{eqnarray}
then for any $\theta\in\Theta^+$, the estimator $\widetilde{\theta}_{\rho_1}$  satisfies
\begin{eqnarray}\label{majoration_risque}
\mathbb{E}_{\theta}\left[l\left(\widetilde{\theta}_{\rho_1},\theta\right)\right] \leq
L_1(K)\inf_{m\in \mathcal{M}}\left[l(\theta_{m,\rho_1},\theta) + \pen(m)\right] +
L_2(K )\frac{\rho_1^2\varphi_{\text{max}}(\Sigma)}{np^2}\ ,
\end{eqnarray}
A similar bound holds if one replaces $\widetilde{\theta}_{\rho_1}$ by  $\widetilde{\theta}_{\rho_1}^{\text{\emph{iso}}}$, $\Theta^+$ by $\Theta^{+,\text{iso}}$, $\theta_{m,\rho_1}$ by $\theta_{m}^{\text{\emph{iso}}}$, and $d_m$ by $d_m^{\text{\emph{iso}}}$.
\end{thrm}
The proof is postponed to Section \ref{section_preuve_mainthrm}. It is based on a novel concentration inequality for suprema of Gaussian chaos stated  in Section \ref{preuves_concentration}.
The constant $K_0$ is made explicit in the proof. Observe that the theorem holds for any $n$, any $p$ and that we have not performed any assumption on the target $\theta\in\Theta^+$ (resp. $\Theta^{+,\text{iso}}$).  If the collection $\mathcal{M}$ does not contain the empty model, one gets the more readable upper bound 
$$\mathbb{E}_{\theta}\left[l\left(\widetilde{\theta}_{\rho_1},\theta\right)\right]\leq
L(K)\inf_{m\in \mathcal{M}}\left[l(\theta_{m,\rho_1},\theta) + \pen(m)\right] \ . $$
This theorem tells us that $\widetilde{\theta}_{\rho_1}$ essentially performs as well as the best trade-off between the bias term $l(\theta_{m,\rho_1},\theta)$ and  $\rho_1^2\varphi_{\text{max}}(\Sigma)\frac{d_m}{np^2}$ that plays the role of a variance. Here are some additional comments. 	~\\~\\
{\bf Remark 1}.  Consider the special case where the target $\theta$ belongs to some parametric set $\Theta^+_m$ with  $m\in\mathcal{M}$. Suppose that  the hypothesis $(\mathbb{H}_1)$ defined in (\ref{def_h1}) is fulfilled.
Choosing a penalty 
$\pen(m) = K\rho_1^2\varphi_{\text{max}}(\Sigma)\frac{d_m}{np^2}$, we get 
\begin{eqnarray}
\mathbb{E}_{\theta}\left[l\left(\widetilde{\theta}_{\rho_1},\theta\right)\right]\leq L(K)\rho_1^2\varphi_{\text{max}}(\Sigma)\frac{d_m}{np^2}\ .\label{risque_parametrique}
\end{eqnarray}
We shall prove in Section \ref{section_risque_asymptotique} and \ref{section_adaptive_sparsity} that this rate is optimal both from an asymptotic oracle  and a minimax point of view. We have mentioned in Section \ref{section_cls} that  $(\mathbb{H}_1)$ is necessary for the bound (\ref{risque_parametrique}) to hold. If $\rho_1$ is chosen large enough, then Assumption $(\mathbb{H}_1)$ is fulfilled. We do not have access to this minimal $\rho_1$ that ensures $(\mathbb{H}_1)$, since it requires the knowledge of $\theta$. Nevertheless, we argue in Section \ref{section_discussion_hypotheses} that ``moderate'' values for $\rho_1$ ensure Assumption $(\mathbb{H}_1)$ when the model $m$ is small.~\\~\\
{\bf Remark 2}. We have mentioned in the introduction that our objective was to obtain oracle inequalities of the form
$$\mathbb{E}_{\theta}\left[l\left(\widetilde{\theta}_{\rho_1},\theta\right)\right]\leq L(K) \inf_{m\in \mathcal{M}}\mathbb{E}\left[l\left(\widehat{\theta}_{m,\rho_1},\theta\right)\right]=L(K)\mathbb{E}\left[\left(\widehat{\theta}_{m^*,\rho_1},\theta\right)\right] \ .$$
This is why we want to compare the sum $l(\theta_{m,\rho_1},\theta)+\pen(m)$ with $\mathbb{E}[l(\widehat{\theta}_{m,\rho_1},\theta)]$.
First, we provide in Section \ref{section_biais} a sufficient condition so that the risk $\mathbb{E}[l(\widehat{\theta}_{m,\rho_1},\theta)]$ decomposes exactly as the sum $l(\theta_{m,\rho_1},\theta)+\mathbb{E}[l(\widehat{\theta}_{m,\rho_1},\theta_{m,\rho_1})]$.
Moreover, we compute in Section \ref{section_risque_asymptotique} the asymptotic variance term $\mathbb{E}[l(\widehat{\theta}_{m,\rho_1},\theta_{m,\rho_1})]$ and compare it with the penalty term $\rho_1^2\varphi_{\text{max}}(\Sigma)\frac{d_m}{np^2}$. We shall then derive oracle type inequalities and discuss the dependency of the different bounds on $\varphi_{\text{max}}(\Sigma)$.~\\~\\
{\bf Remark 3}. Condition (\ref{condition_penalite}) gives a lower bound on the penalty function $\pen(.)$ so that the result holds. Choosing a proper penalty term according to (\ref{condition_penalite}) therefore requires an upper bound on the largest eigenvalue of $\Sigma$. However, such a bound is seldom known in practice. We shall mention in Section \ref{section_discussion} a practical method to calibrate the penalty.\\

A  bound similar to (\ref{majoration_risque}) holds for the Frobenius distance between the partial correlation matrices   ($I_{p^2}-C(\theta)$) and ($I_{p^2}-C(\widetilde{\theta}_{\rho_1})$).
\begin{cor}\label{corollaire_Frobenius}
Assume the same as in Theorem \ref{mainthrm}, except that there is equality in (\ref{condition_penalite}). Then,
\begin{eqnarray}
\mathbb{E}_{\theta}\left[\|C(\widetilde{\theta}_{\rho_1})-C\left(\theta\right)\|_F^2\right] &\leq&L_1\left(K\right)
\frac{\varphi_{\text{\emph{max}}}(\Sigma)}{\varphi_{\text{\emph{min}}}(\Sigma)}\inf_{m\in \mathcal{M}}\left[\|C(\theta_{m,\rho_1})-C(\theta)\|_F^2 + \frac{K\rho_1^2d_m}{n}\right]\nonumber\\  &+ &L_2(K)\frac{\varphi_{\text{\emph{max}}}(\Sigma)}{\varphi_{\text{\emph{min}}}(\Sigma)}\frac{\rho_1^2}{n}\ .\label{majoration_risque_frobenius}
\end{eqnarray}
A similar result holds for isotropic GMRFs.
\end{cor}

\begin{proof}[Proof of Corollary \ref{corollaire_Frobenius}] 
This is a consequence of Theorem \ref{mainthrm}.
 By definition (\ref{def_lmatriciel}) of the loss function $l(.,.)$, the two following  bounds hold 
\begin{eqnarray*}
p^2l(\theta_1,\theta_2)\geq \varphi_\text{min}(\Sigma)\|C(\theta_1)-C(\theta_2)\|_F^2\  \\ 
p^2l(\theta_1,\theta_2)\leq \varphi_\text{max}(\Sigma)\|C(\theta_1)-C(\theta_2)\|_F^2\ .
\end{eqnarray*}
Gathering these bounds with (\ref{majoration_risque}) yields the result.
\end{proof}

 The same comments as for Theorem (\ref{mainthrm}) hold. We may express this Corollary \ref{corollaire_Frobenius} in terms of the risk $\mathbb{E}(\|\widetilde{\theta}_{\rho_1}-\theta\|_F^2)$, since $\|C(\theta_1)-C(\theta_2)\|_F^2=p^2\|\theta_1-\theta_2\|_F^2$: 
\begin{eqnarray*}
\mathbb{E}_{\theta}\left[\|\widetilde{\theta}_{\rho_1}-\theta\|_F^2\right] &\leq&L_1\left(K\right)
\frac{\varphi_{\text{\emph{max}}}(\Sigma)}{\varphi_{\text{\emph{min}}}(\Sigma)}\inf_{m\in \mathcal{M}}\left[\|\theta_{m,\rho_1}-\theta\|_F^2 + \frac{K\rho_1^2d_m}{np^2}\right]\nonumber\\  &+ &L_2(K)\frac{\varphi_{\text{\emph{max}}}(\Sigma)}{\varphi_{\text{\emph{min}}}(\Sigma)}\frac{\rho_1^2}{np^2}\ .
\end{eqnarray*}

\section{Parametric risk and asymptotic oracle inequalities}\label{section_risk}

In this section, we study the risk of the parametric estimators $\widehat{\theta}_{m,\rho_1}$ in order to assess the optimality of Theorem \ref{mainthrm}. 

\subsection{Bias-variance decomposition}\label{section_biais}

The properties of the parametric estimator $\widehat{\theta}_{m,\rho_1}$ and of the projection $\theta_{m,\rho_1}$ differ slightly whether $\theta_{m,\rho_1}$ belongs to the open set $\Theta_{m,\rho_1}^{+}$ or to its border. Observe that Hypothesis $(\mathbb{H}_1)$ defined in (\ref{def_h1}) does not necessarily imply that the projection $\theta_{m,\rho_1}$ belongs to $\Theta_{m}^{+}$. This is why we introduce the condition $(\mathbb{H}_2)$:
\begin{eqnarray}\label{def_h2}
\theta\in\mathcal{B}_1(0_p,1)\hspace{1cm}\Longleftrightarrow\hspace{1cm} \|\theta\|_1<1\  .  
\end{eqnarray}
The condition $\|\theta\|_1<1$ is equivalent to $(I_{p^2}-C(\theta))$ is strictly diagonally dominant. 
Condition  $(\mathbb{H}_2)$ implies that the largest eigenvalue of $(I_{p^2}-C(\theta))$ is smaller than 2 and therefore that  $(\mathbb{H}_1)$ is fulfilled since $\rho_1$ is supposed larger than $2$.
We further discuss this assumption in Section \ref{section_discussion_hypotheses}.

\begin{lemma}\label{lemme_projection}
Let $\theta\in \Theta^+$ such that $(\mathbb{H}_2)$ holds and let $m\in \mathcal{M}_1$. Then, the minimum of  $\gamma(.)$ over $\Theta_m$ is achieved in $\Theta^+_{m,2}$. This implies that 
$$\theta_{m,\rho_1} = \arg\min_{\theta'\in \Theta_m} \gamma(\theta')\, \hspace{1.5cm}\text{and}\hspace{1.5cm}\gamma(\theta_{m,\rho_1})=\var_{\theta}\left(X{\scriptstyle[0,0]}|X_m\right)\ .$$
Besides, $\|\theta_{m,\rho_1}\|_1\leq \|\theta\|_1$. The same results holds
for  $\theta_{m,\rho_1}^{\text{iso}}$ if $\theta$ in $\Theta^{+,\text{\emph{iso}}}$.
\end{lemma}
The proof is given in the technical appendix~\cite{technical}. 
The purpose of this property is threefold. First, we derive that Assumption $(\mathbb{H}_2)$ ensures that $\theta_{m,\rho_1}$ belongs $\Theta_{m,\rho_1}^+$ and that the smallest eigenvalue of $(I_{p^2}-C(\theta_{m,\rho_1}))$ is larger than $1-\|\theta\|_1$.
Second, it allows to express the
projection $\theta_{m,\rho_1}$ in terms of conditional expectation (Corollary
\ref{definition_thetam}). Finally,  we deduce 
a bias-variance decomposition of the estimator $\widehat{\theta}_{m,\rho_1}$ (Corollary \ref{decomposition_biais-variance}).
In other words, the equality holds in (\ref{decomposition_biais_variance_inegalite}).

\begin{cor}\label{definition_thetam}
Let $\theta\in \Theta^+$ such that $(\mathbb{H}_2)$ holds and let $m\in \mathcal{M}_1$. The projection  $\theta_{m,\rho_1}$ is
uniquely defined by the equation
\begin{eqnarray*}
\mathbb{E}_{\theta}\left(X{\scriptstyle[0,0]}|X_{ m }\right)= \sum_{(i,j)\in
   m }\theta_{m,\rho_1}{\scriptstyle[i,j]}X{\scriptstyle[i,j]}\ ,
\end{eqnarray*}
and $\theta_{m,\rho_1}{\scriptstyle[i,j]}=0$ for any $(i,j)\notin  m $. Similarly, if $\theta\in
\Theta^{+,\text{\emph{iso}}}$ satisfies $(\mathbb{H}_2)$, then $\theta^{\text{\emph{iso}}}_{m,\rho_1}$ is uniquely
defined by the equation
\begin{eqnarray*}
\mathbb{E}_{\theta}\left(X{\scriptstyle[0,0]}|X_{ m }\right)= \sum_{(i,j)\in
  m }\theta^{\text{\emph{iso}}}_{m,\rho_1}{\scriptstyle[i,j]}X{\scriptstyle[i,j]}\ ,
\end{eqnarray*}
and $\theta^{\text{\emph{iso}}}_{m,\rho_1}{\scriptstyle[i,j]}=0$ for any $(i,j)\notin  m $.
\end{cor}

Consequently,
$\sum_{1\leq i,j\leq p} \theta_{m,\rho_1}{\scriptstyle[i,j]}X{\scriptstyle[i,j]}$ is the best linear predictor of
$X{\scriptstyle[0,0]}$ given the covariates $X{\scriptstyle[i,j]}$ with $(i,j)\in m$. This is precisely the definition of the kriging parameters (Stein~\cite{stein}). Hence, the matrix $\theta_{m,\rho_1}$ corresponds to the kriging parameters of $X{\scriptstyle[0,0]}$ with kriging  neighborhood's range  of $r_m$. The distance $r_m$ is introduced  in Definition \ref{definition_modele} and stands for the radius of $m$.

\begin{cor}\label{decomposition_biais-variance}
Let $\theta\in \Theta^+$  such that $(\mathbb{H}_2)$ holds  and let $m\in\mathcal{M}_1$. 
The loss of $\widehat{\theta}_{m,\rho_1}$ decomposes as
 $l(\widehat{\theta}_{m,\rho_1},\theta) =
l(\theta_{m,\rho_1},\theta)+ l(\widehat{\theta}_{m,\rho_1},\theta_{m,\rho_1})$. If $\theta$ belongs to $\Theta_m^{+,\text{\emph{iso}}}$ and $(\mathbb{H}_2)$ holds, then we also have the decomposition 
 $l(\widehat{\theta}^{\text{\emph{iso}}}_{m,\rho_1},\theta) =
l(\theta^{\text{\emph{iso}}}_{m,\rho_1},\theta)+ l(\widehat{\theta}^{\text{\emph{iso}}}_{m,\rho_1},\theta_{m,\rho_1})$.
\end{cor}
A proof is provided in the technical appendix~\cite{technical}.
If $\theta$ does not satisfy Assumption $(\mathbb{H}_2)$, then $\theta_{m,\rho_1}$ does not necessarily belong to $\Theta_{m,\rho_1}^{+}$ and there may not be such a bias variance decomposition.

\subsection{Asymptotic risk}\label{section_risque_asymptotique}

In this section, we evaluate the risk of each estimator $\widehat{\theta}_{m,\rho_1}$ and use it as a benchmark to assess the result of Theorem \ref{mainthrm}. We have mentioned in Corollary \ref{decomposition_biais-variance} that under $(\mathbb{H}_2)$ the risk $\mathbb{E}_{\theta}[l(\widehat{\theta}_{m,\rho_1},\theta)]$ decomposes into the sum of the bias $l(\theta_{m,\rho_1},\theta)$ and a variance term $\mathbb{E}_{\theta}[l(\widehat{\theta}_{m,\rho_1},\theta_{m,\rho_1})]$. If this last quantity is of the same order as the penalty $\pen(m)$ introduced in (\ref{condition_penalite}), then Theorem \ref{mainthrm} yields an  oracle inequality.
However, we are unable to express this variance term $\mathbb{E}_{\theta}[l(\widehat{\theta}_{m,\rho_1},\theta_{m,\rho_1})]$ in a simple form. This is why we restrict ourselves to study the risks when  $n$ tends to infinity. Nevertheless, these results give us some hints to appreciate the strength and the weaknesses of Theorem \ref{mainthrm} and the upper bound (\ref{risque_parametrique}).

In the following proposition, we adapt a result of Guyon \citep{guyon95} Sect.4.3.2 to obtain an asymptotic expression of the risk $\mathbb{E}_{\theta}[l(\widehat{\theta}_{m,\rho_1},\theta_{m,\rho_1})]$. We first need to introduce some new notations. For any model $m$ in the collection $\mathcal{M}_1\setminus\{\emptyset\}$, we fix a sequence $(i_k,j_k)_{k=1\ldots d_m}$ of integers such that $(\Psi_{i_1,j_1},\ldots, \Psi_{i_{d_m},j_{d_m}})$ is a basis of the space $\Theta_m$. Then,  $\chi_m{\scriptstyle[0,0]}$ stands for the random vector of size $d_m$ that contains the neighbors of $X{\scriptstyle[0,0]}$ 
$$\chi_m{\scriptstyle[0,0]}^* := \left[tr\left(\Psi_{i_1,j_1}X^v\right),\ldots, tr\left(\Psi_{i_{d_m},j_{d_m}}X^v\right)\right]\ .$$
Besides, for any $\theta\in \Theta^+$, we define the matrices $V$, $W$ and \emph{IL}$_m$ as
\begin{eqnarray*}\left\{\begin{array}{ccc}
 V & := &  \cov_{\theta}(\chi_m{\scriptstyle[0,0]})\\
W{\scriptstyle[k,l]} & := & \frac{1}{p^2}tr\left[C\left(\Psi_{i_k,j_k}\right)\left(I_{p^2}-C(\theta_{m,\rho_1})\right)^2 \left(I_{p^2}-C(\theta)\right)^{-2}C\left(\Psi_{i_l,j_l}\right)\right],\text{ for any } k=1,\ldots,d_m\ \\
 \text{\emph{IL}}_m  & := & \text{Diag}\left(\|\Psi_{i_k,j_k}\|_F^2,\ k=1,\ldots,d_m\right) \ ,\end{array}\right.
\end{eqnarray*}	
where  for any vector $u$, $\text{Diag}(u)$ is the diagonal matrix whose diagonal elements are the components of $u$. We also define the corresponding quantities $\chi^{\text{iso}}_m{\scriptstyle[0,0]}$, $V^{\text{iso}}$, $W^{\text{iso}}$, and $\text{\emph{IL}}^{\text{iso}}_m$ in order to consider the isotropic estimator $\widehat{\theta}_{m,\rho_1}^{\text{iso}}$.

\begin{prte}\label{variance_vrai}
Let  $m$ be a model in  $\mathcal{M}_1\setminus\{\emptyset\}$ and let $\theta$ be an element of $\Theta^{+}_{m}$ that satisfies $(\mathbb{H}_1)$. Then,  $\widehat{\theta}_{m,\rho_1}$ converges to $\theta$ in probability and
\begin{eqnarray}
\lim_{n\rightarrow +\infty}np^2 \mathbb{E}_{\theta}\left[l\left(\widehat{\theta}_{m,\rho_1},\theta\right)\right]
= 2\sigma^4tr\left[\text{IL}_mV^{-1}\right]\ . \label{risqueas2}
\end{eqnarray}
Let $\theta$ in $\Theta^+$ such that  $(\mathbb{H}_2)$ is fulfilled. Then,  $\widehat{\theta}_{m,\rho_1}$ converges to $\theta_{m,\rho_1}$ in probability and 
\begin{eqnarray}	
\lim_{n\rightarrow +\infty}np^2
\mathbb{E}_\theta\left[l\left(\widehat{\theta}_{m,\rho_1},\theta_{m,\rho_1}\right)\right]= 2\sigma^4tr(WV^{-1})\ . \label{risqueas1}
\end{eqnarray}
Both results still hold for the estimator $\widehat{\theta}_{m,\rho_1}^{\text{\emph{iso}}}$ if  $\theta$ belongs to $\Theta^{+,\text{\emph{iso}}}$ and if one replaces $V$, $W$, and IL$_m$  by $V^{\text{\emph{iso}}}$, $W^{\text{\emph{iso}}}$,  and IL$_m^{\text{\emph{iso}}}$.
\end{prte}
In the first case, Assumption $(\mathbb{H}_1)$ ensures that $\theta\in\Theta^+_{m,\rho_1}$, whereas Assumption $(\mathbb{H}_2)$ ensures that $\theta_{m,\rho_1}\in\Theta^+_{m,\rho_1}$. The proof is based on the extension of Guyon's approach in the toroidal framework.

The expressions (\ref{risqueas2}) and (\ref{risqueas1}) are not easily interpretable in the present form. This is why we first derive (\ref{risqueas2}) when $\theta$ is zero. Observe that it is equivalent to the independence of the $(X{\scriptstyle[i,j]})_{(i,j)\in \Lambda}$.
\begin{ex}\label{variance_identite}
Assume that $\theta$ is zero. Then, for any model
$m\in \mathcal{M}_1$, the asymptotic risks of $\widehat{\theta}_{m,\rho_1}$ and $\widehat{\theta}^{\text{\emph{iso}}}_{m,\rho_1}$ satisfy
\begin{eqnarray*}
\lim_{n\rightarrow +\infty}np^2
\mathbb{E}_{0_p}\left[l\left(\widehat{\theta}_{m,\rho_1},0_p\right)\right]= 2\sigma^2d_m\vspace{0.5cm}
\text{ and }\vspace{0.5cm}\lim_{n\rightarrow +\infty}np^2
\mathbb{E}_{0_p}\left[l\left(\widehat{\theta}^{\text{\emph{iso}}}_{m,\rho_1},0_p\right)\right]= 2\sigma^2d^{\text{\emph{iso}}}_m\ ,
\end{eqnarray*} 
where we recall that $d^{\text{\emph{iso}}}_m$  is the dimension of the space $\Theta_{m}^{\text{\emph{iso}}}$.
\end{ex}
\begin{proof}
 Since the components of $X$ are independent, the matrix $V$ equals $\sigma^2 \text{\emph{IL}}_m$. We conclude by applying Proposition \ref{variance_vrai}
\end{proof}
Therefore, when the variables $X{\scriptstyle[i,j]}$ are independent, the asymptotic risk of $\widehat{\theta}_{m,\rho_1}$ equals, up to a factor 2,  the variance term of the least squares estimator in the fixed design Gaussian regression framework. This quantity is of the same order as the penalty introduced in Section \ref{section_main_result}. When the matrix $\theta$ is non zero, we  can lower bound the limits (\ref{risqueas2}) and (\ref{risqueas1}).

\begin{cor}\label{minoration_variance_asymptotique}
Let $m$ be a model in $\mathcal{M}_1$ and let $\theta\in \Theta_m^+$ that satisfies $(\mathbb{H}_1)$.    
Then, the variance term is asymptotically lower bounded as follows
\begin{eqnarray}\label{minorationas}
\lim_{n\rightarrow +\infty}np^2
\mathbb{E}_\theta\left[l\left(\widehat{\theta}_{m,\rho_1},\theta\right)\right]\geq  L\sigma^2\varphi_{\text{min}}\left[{I_{p^2}}-C(\theta)\right]d_m=L\sigma^4\frac{d_m}{\varphi_{\text{max}}(\Sigma)}
\ ,
\end{eqnarray} 
where $L$ is a universal constant. Let $\theta\in \Theta^+$ that satisfies $(\mathbb{H}_2)$. For any model $m\in\mathcal{M}_1$, 
\begin{eqnarray}\label{minorationas2}
\lim_{n\rightarrow +\infty}np^2
\mathbb{E}_\theta\left[l\left(\widehat{\theta}_{m,\rho_1},\theta_{m,\rho_1}\right)\right]\geq L\sigma^2\left(1-\|\theta\|_1\right)^3d_m\ ,
\end{eqnarray} 

\end{cor}
The proof is postponed to the technical appendix~\cite{technical}. Again, analogous lower bounds  hold for $\widehat{\theta}^{\text{iso}}_{m,\rho_1}$ when $\theta$ belongs to $\Theta^{\text{iso},+}$. 
This corollary states that asymptotically with respect to $n$ the variance term of $\widehat{\theta}_{m,\rho_1}$ is larger than the order $d_m/(np^2)$. This expression is not really surprising since $d_m$ stands for the dimension of the model $m$ and $np^2$ corresponds to the number of data observed.
Let define $R_{\theta,\infty}(\widehat{\theta}_{m,\rho_1},\theta_{m,\rho_1}):= \lim_{n\rightarrow +\infty}np^2\mathbb{E}_{\theta}[l(\widehat{\theta}_{m,\rho_1},\theta_{m,\rho_1})]$ as the asymptotic variance term for $\widehat{\theta}_{m,\rho_1}$ rescaled by the number $np^2$ of observations.\\

The first part of the corollary (\ref{minorationas}) states that from an asymptotic point of view the upper bound (\ref{risque_parametrique}) is optimal. By Theorem \ref{mainthrm}, if we choose  $\pen(m)=K\rho_1^2\varphi_{\text{max}}(\Sigma)\frac{d_m}{np^2}$, then it holds that
$$\mathbb{E}\left[l\left( \widetilde{\theta}_{\rho_1},\theta \right)\right]\leq L\left(K,\rho_1,\varphi_{\text{min}}\left[I_{p^2}-C(\theta)\right]\right)\frac{R_{\theta,\infty}(\widehat{\theta}_{m,\rho_1},\theta)}{np^2} ,$$
for any model $m\in\mathcal{M}\setminus\emptyset$ and any $\theta\in\Theta_m^+$ that satisfies $(\mathbb{H}_1)$. This property holds for any $n$ and any $p$. Hence, $\widetilde{\theta}_{\rho_1}$ performs as well as the parametric estimator $\widehat{\theta}_{m,\rho_1}$ if the support of $\theta$ belongs to  some unknown model $m$ and if $\theta$ satisfies $(\mathbb{H}_1)$.\\

If we assume that $\|\theta\|_1< 1$ (Hypothesis  $(\mathbb{H}_2)$), we are able to derive a stronger result.
\begin{prte}\label{oracle_asymptotique}
Considering $K\geq K_0$, $\rho_1\geq 2$, $\eta<1$ and a collection $\mathcal{M}\subset\mathcal{M}_1\setminus\emptyset$, we define the estimator $\widetilde{\theta}_{\rho_1}$ with the penalty $\pen(m)=K\rho_1^2\frac{d_m}{np^2(1-\eta)}$. Then, the risk of $\widetilde{\theta}_{\rho_1}$ is upper bounded by  
\begin{eqnarray}\label{inegalite_oracle_asymptotique}
\mathbb{E}_{\theta}\left[l\left(\widetilde{\theta}_{\rho_1},\theta\right)\right]\leq L(K,\rho_1,\eta) \inf_{m\in\mathcal{M}}\left\{l\left(\theta_{m,\rho_1},\theta\right) + \frac{R_{\theta,\infty}\left(\widehat{\theta}_{m,\rho_1},\theta_{m,\rho_1}\right)}{np^2}\right\} \ ,
\end{eqnarray}
for any $\theta\in \Theta^+\cap \mathcal{B}_1\left(0_p,\eta\right)$.
\end{prte}
Observe that this property holds for any $n$ and any $p$. If the matrix $\theta$ is strictly diagonally dominant, we therefore obtain an upper bound similar to an oracle inequality, except that the variance term $\mathbb{E}_{\theta}[l(\widehat{\theta}_{m,\rho_1},\theta_{m,\rho_1})]$ has been replaced by its asymptotic counterpart $R_{\theta,\infty}(\widehat{\theta}_{m,\rho_1},\theta_{m,\rho_1})/(np^2)$. However, this  inequality is not valid uniformly over any $\eta<1$ : when $\eta$ converges to one, the constant $L(K,\rho_1,\eta)$ tends to infinity. Indeed, if $\|\theta\|_1$ converges to one, the lower bound (\ref{minorationas2}) on the variance term  can behave like $(1-\|\theta\|_1)^3d_m/(np^2)$ for some matrices $\theta$ whereas the penalty term $d_m/[np^2(1-\|\theta\|_1)]$ tends to infinity. \\

In the remaining part of the section, we illustrate that the constant $L(K,\eta,\rho_1)$ has to go to infinity when $\eta$ goes to one. Let us consider the  model $m_1$. It consists of GMRFs with 4-nearest neighbors.
\begin{ex}\label{variance_limite_alpha}
Let  $\theta$ be a non zero element of $\Theta_{m_1}^{\text{\emph{iso}}}$, then the asymptotic risk of $\widehat{\theta}_{m_1,\rho_1}^{\text{\emph{iso}}}$ simplifies as
\begin{eqnarray}\label{variance_simple_m1}
\lim_{n\rightarrow
  +\infty}np^2\mathbb{E}_{\theta}\left[l\left(\widehat{\theta}_{m_1,\rho_1}^{\text{\emph{iso}}},\theta \right)\right]= 2 \frac{\sigma^4\theta{\scriptstyle[1,0]}}{\cov(X{\scriptstyle[1,0]},X{\scriptstyle[0,0]})}\ .
\end{eqnarray}
If we let the size $p$ of the network tend to infinity and $\theta{\scriptstyle[1,0]}$ go to $1/4$, the risk is equivalent to
\begin{eqnarray*}
\lim_{p\rightarrow  +\infty} \lim_{n\rightarrow
  +\infty}np^2\mathbb{E}_{\theta}\left[l\left(\widehat{\theta}_{m_1,\rho_1}^{\text{\emph{iso}}},\theta \right)\right]\begin{array}{c}\sim\\{\tiny \theta{\scriptstyle[1,0]}\rightarrow 1/4}\end{array}
 \frac{16\sigma^2(1-4\theta{\scriptstyle[1,0]})}{\log(16)} \ .
\end{eqnarray*}
\end{ex}
The proof is postponed to the technical appendix~\cite{technical}.
If follows from the second result that the lower bound (\ref{minorationas}) is sharp since in this particular case $\varphi_{\text{min}}(I_{p^2}-C(\theta))=\sigma^2(1-4\theta{\scriptstyle[1,0]})$.  When $\theta{\scriptstyle[1,0]}$ tends to $1/4$, then $\|\theta\|_1$ tends to one and $\mathbb{E}_{\theta}[l(\widehat{\theta}_{m_1,\rho_1}^{\text{\emph{iso}}},\theta )]$ behaves like $\sigma^2(1-\|\theta\|_1)d_{m_1}^{\text{iso}}/(np^2)$ whereas the penalty $\pen(m_1)$ given in Theorem \ref{mainthrm} has to be larger than $\sigma^2d_{m_1}^{\text{iso}}/[np^2(1-\|\theta\|_1)]$. Hence, the variance term and the penalty $\pen(.)$ are not necessarily of the same order when $\|\theta\|_1$ tends to one.
Theorem \ref{mainthrm} cannot lead to an oracle inequality of the type (\ref{inegalite_oracle_asymptotique}), which is valid uniformly on $\eta<1$.

\begin{ex}\label{variance_exemple_risque_gros}
Let $\alpha$ be a positive number smaller than $1/4$. For any integer $p$ which is divisible by 4, we define the $p\times p$ matrix $\theta^{(p)}$ by
$$\left\{\begin{array}{ccc}\theta^{(p)}{\scriptstyle[p/4,p/4]} =\theta^{(p)}{\scriptstyle[-p/4,p/4]}=\theta^{(p)}{\scriptstyle[p/4,-p/4]}=\theta^{(p)}{\scriptstyle[-p/4,-p/4]} & := & \alpha \\ \theta^{(p)}{\scriptstyle[i,j]}& := & 0\, \, \text{\ else}\  .\end{array}\right.$$ 
Then, the variance term is asymptotically lower bounded as follows
\begin{eqnarray*}
 \lim_{p\rightarrow  +\infty} \lim_{n\rightarrow
  +\infty}np^2\mathbb{E}_{\theta^{(p)}}\left[l\left(\widehat{\theta^{(p)}}_{m_1,\rho_1}^{\text{\emph{iso}}},[\theta^{(p)}]_{m_1,\rho_1}^{\text{\emph{iso}}} \right)\right]\geq \frac{L\sigma^2}{1-4\alpha} \ .
\end{eqnarray*}
\end{ex}
The proof is postponed to the technical appendix~\cite{technical}. This variance term is of order $\sigma^2d^{\text{iso}}_m/[np^2(1-\|\theta\|_1)]=\varphi_{\text{max}}(\Sigma)d^{\text{iso}}_m/(np^2)$ when $\|\theta\|_1$ goes to one. 
The penalty $\pen(m)$ introduced in Proposition \ref{oracle_asymptotique} is therefore a sharp upper bound of the variance terms.

On the one hand, we take a penalty $\pen(m)$ larger than $\sigma^2d_m/(np^2(1-\|\theta\|_1))$. On the other hand, the variance of $\widehat{\theta}_{m,\rho_1}$ is of  the order $\sigma^2(1-\|\theta\|_1)d_m/(np^2)$ in some cases. The bound (\ref{inegalite_oracle_asymptotique}) cannot therefore hold uniformly over any $\eta<1$. We think that it is intrinsic to the penalization strategy.

\section{Comments on the assumptions}\label{section_discussion_hypotheses}

%%%%%%%%%%%%%%%%%%%%%%%%%%%%%%%%%%%%%%%%%%%
%%%%  Discussion de rho1 et thetha
%%%%%%%%%%%%%%%%%%%%%%%%%%%%%%%%%%%%%%%%%%

In this section, we discuss the dependency of the estimators $\widehat{\theta}_{m,\rho_1}$ on $\rho_1$ as well as Assumptions $(\mathbb{H}_1)$ and $(\mathbb{H}_2)$.~\\~\\	
\textbf{Dependency of $\widehat{\theta}_{m,\rho_1}$ on $\rho_1$}. We recall that the estimator $\widehat{\theta}_{m,\rho_1}$ is defined in (\ref{definition_estimateur_cls}) as the minimizer of the CLS empirical contrast $\gamma_{n,p}(.)$ over $\Theta^+_{m,\rho_1}$. 
It may seem restrictive to perform the minimization over the set $\Theta^+_{m,\rho_1}$ instead of $\Theta^+_{m}$. Nevertheless,  we advocate that it is not the case, at least for small models. Let us indeed define  
\begin{eqnarray*}
 \rho(m):=\sup_{\theta\in \Theta_m^+}\varphi_{\text{max}}\left[I_{p^2}-C(\theta)\right]\text{\hspace{0.5cm} and \hspace{0.5cm}}\rho^{\text{\emph{iso}}}(m):=\sup_{\theta\in \Theta_m^{+,\text{\emph{iso}}}}\varphi_{\text{max}}\left[I_{p^2}-C(\theta)\right]\ .
\end{eqnarray*}
The quantities $\rho(m)$ and $\rho^{\text{\emph{iso}}}(m)$ are finite since $\Theta_m^+$ is bounded. If one takes $\rho_1$ larger than $\rho(m)$ (resp. $\rho^{\text{iso}}(m)$), then the set $\Theta_{m,\rho_1}^+$ (resp. $\Theta_{m,\rho_1}^{+,\text{iso}}$) is exactly $\Theta_{m}^+$ (resp. $\Theta_{m}^{+,\text{iso}}$). We illustrate in Table \ref{tableau_comparaison_rho} that $\rho(m)$ and $\rho^{\text{iso}}(m)$ are small, when the model $m$ is small. Consequently, choosing a moderate value for $\rho_1$ is not really restrictive for small models. However, when the size of the model $m$ increases, the sets $\Theta_{m,\rho_1}^{+}$ and $\Theta_{m}^{+}$ become different for moderate values of $\rho_1$. In Section \ref{section_discussion}, we discuss the choice of $\rho_1$.
\begin{Table}[h]
\caption{Approximate computation of $\rho(m)$ and $\rho^{\text{iso}}(m)$ for the four smallest models with $p=50$.\label{tableau_comparaison_rho}}
\begin{center}
\begin{tabular}{|l|c|c|c|c|}
  \hline $d_m$   & 2 & 4 & 6 & 10 \\ \hline
      $\rho(m)$ & $2.0$ & $4.0$ & $5.0$ & $6.8$   \\ \hline  
     $d_{m}^{\text{iso}}$& 1 & 2&3 &4  \\ \hline 
     $\rho^{\text{iso}}(m)$  &  $2.0$ &  $ 4.0$ & $5.0$& $6.8$ \\ \hline 
\end{tabular}
\end{center}
\end{Table}
%%%%%%%%%%%%%%%%%%%%%%%%%%%%%%%%%%%%%%%%%%%%%%%%%%%%%%%%%%
% Discussion H1
%%%%%%%%%%%%%%%%%%%%%%%%%%%%%%%%%%%%%%%%%%%%%%%%%%%%%%%%%%%
~\\ ~\\ \textbf{Assumption $(\mathbb{H}_1)$} defined in (\ref{def_h1}) states that the largest eigenvalue of $(I_{p^2}-C(\theta))$ is smaller than $\rho_1$. We have illustrated in Table \ref{tableau_comparaison_rho} that if the support of $\theta$ belongs to a small model $m$, then the maximal absolute value of $(I_{p^2}-C(\theta))$ is  small. Hence, Assumption $(\mathbb{H}_1)$ is ensured for ``moderate'' values of $\rho_1$ as soon as the support of $\theta$ belongs to some small model. 
If $\theta$ is not sparse but approximately sparse it is likely that the largest eigenvalue of $\theta$ remain moderate. In practice, we do not know in advance if a given choice of $\rho_1$ ensures $(\mathbb{H}_1)$. In Section \ref{section_discussion}, we discuss an extension of our procedure which does not require Assumption $(\mathbb{H}_1)$.~\\~\\
%%%%%%%%%%%%%%%%%%%%%%%%%%%%%%%%%%%%%%%%%%%%%%%%%%%%%%%%%%
% Discussion H2
%%%%%%%%%%%%%%%%%%%%%%%%%%%%%%%%%%%%%%%%%%%%%%%%%%%%%%%%%%%
\textbf{Assumption $(\mathbb{H}_2)$} defined in (\ref{def_h2}) states that $\theta\in \mathcal{B}_1(0_p,1)$ or equivalently that the matrix $(I_{p^2}-C(\theta))$ is diagonally dominant. Rue and Held prove in \cite{rue} Sect.2.7 that $\Theta_{m_1}^+$ is included in $\mathcal{B}_1(0_p,1)$. They also point out that a small part of $\Theta_{m_2}^+$ does not belong to $\mathcal{B}_1(0_p,1)$. In fact, Assumption $(\mathbb{H}_2)$ becomes more and more restrictive if the support of $\theta$ becomes larger.  Nevertheless, Assumption $(\mathbb{H}_2)$ is also quite common in the literature (as for instance  in \cite{guyon95}).

If one looks closely at our proofs involving Assumptions $(\mathbb{H}_2)$, one realizes that this assumptions is only made to ensure the following facts:
\begin{enumerate}
\item The \emph{projection} $\theta_{m,\rho_1}$ belongs to the open set $\Theta^+_{m,\rho_1}$ for any model $m\in\mathcal{M}$ (Corollary \ref{decomposition_biais-variance}).
\item The smallest eigenvalue of ($I_{p^2}-C(\theta_{m,\rho_1})$) is  lower bounded by some positive number $\rho_2$,.uniformly over all models $m\in\mathcal{M}$.
\end{enumerate}
From empirical observations, these two last facts seem far more restrictive than $(\mathbb{H}_2)$. We used Assumption  $(\mathbb{H}_2)$ in the statement of our results, because we did not find any weaker but still simple condition that ensures facts 1 and 2.

\section{Minimax rates}\label{section_minimax}

In Theorem \ref{mainthrm} and Proposition \ref{oracle_asymptotique} we have shown that under mild assumptions on $\theta$ the estimator $\widetilde{\theta}_{\rho_1}$ behaves almost as well as the best estimator among the family $\{\widehat{\theta}_{m,\rho_1},\, m\in\mathcal{M}\}$. We now compare the risk of $\widetilde{\theta}_{\rho_1}$ with the risk of any other possible estimator $\widehat{\theta}$. This includes comparison with maximum likelihood methods. There is no hope to make a pointwise comparison with an arbitrary estimator. Therefore, we classically consider the maximal risk over some suitable subsets $\mathcal{T}$ of $\Theta^+$. The \emph{minimax risk} over the set $\mathcal{T}$ is given by $\inf_{\widehat{\theta}}\sup_{\theta\in \mathcal{T}}\mathbb{E}_{\theta}[l(\widehat{\theta},\theta)]$, where the infimum is taken over all possible estimators $\widehat{\theta}$ of $\theta$. Then, the estimator $\widetilde{\theta}_{\rho_1}$ is said to be \emph{approximately minimax} with respect to the set $\mathcal{T}$ if the ratio
\begin{eqnarray*}
 \frac{\sup_{\theta\in \mathcal{T}}\mathbb{E}_{\theta}\left[l\left(\widetilde{\theta}_{\rho_1},\theta\right)\right]}{\inf_{\widehat{\theta}}\sup_{\theta\in \mathcal{T}}\mathbb{E}_{\theta}\left[l\left(\widehat{\theta},\theta\right)\right]}
\end{eqnarray*}
is smaller than a constant that does not depend on $\sigma^2$, $n$ or $p$. An estimator is said to be \emph{adaptive} to a collection $(\mathcal{T}_i)_{i\in\mathcal{I}}$ if it is simultaneously minimax over 
each $\mathcal{T}_i$. The problem of designing  adaptive estimation procedures is in general difficult. It has been  extensively studied in the fixed design Gaussian regression framework. See for instance \cite{birge2001} for a detailed discussion. In the sequel, we  adapt some of their ideas to the GMRF framework.\\

We prove in Section \ref{section_adaptive_sparsity} that the estimator $\widetilde{\theta}_{\rho_1}$ is adaptive to the unknown sparsity of the matrix $\theta$. Moreover, it is also adaptive  if we consider the Frobenius distance between partial correlation matrices. In Section \ref{section_adaptive_ellipsoids}, we show that $\widetilde{\theta}_{\rho_1}$ is also adaptive to the rates of decay of the bias.

We need to restrain ourselves to set of matrices $\theta$ such that the largest eigenvalue of the covariance matrix $\Sigma$ is uniformly bounded. This is why  we define 
\begin{eqnarray}\label{definition_U}
 \forall \rho_2>1\ , \hspace{0.3cm}\mathcal{U}(\rho_2):= \left\{\theta\in\Theta, \varphi_{\text{min}}\left(I_{p^2}-C(\theta)\right)\geq \frac{1}{\rho_2}\right\}\ .
\end{eqnarray}
Observe that $\theta\in \mathcal{U}(\rho_2)$ is exactly equivalent to $\varphi_{\text{max}}(\Sigma)\leq \sigma^2\rho_2$ since $\Sigma=\sigma^2(I_{p^2}-C(\theta))$.

\subsection{Adapting to unknown sparsity}\label{section_adaptive_sparsity}

In this subsection, we prove that under mild assumptions the penalized estimator $\widetilde{\theta}_{\rho_1}$ is adaptive to the unknown sparsity of $\theta$. We first lower bound the minimax rate of convergence on given hypercubes.

\begin{defi}\label{definition_hypercube}
Let $m$ be a  model in the collection $\mathcal{M}_1\setminus\emptyset$. We consider $\left(\Psi_{i_1,j_1},\ldots, \Psi_{i_{d_m},j_{d_m}}\right)$ a basis of the space $\Theta_m$ defined by (\ref{definition_psi}). For any $\theta'\in\Theta_m^+$, the hypercube $\mathcal{C}_m(\theta',r)$ is defined as 
$$\mathcal{C}_m\left(\theta',r\right):= \left\{\theta'+\sum_{k=1}^{d_m}\Psi_{i_k,j_k}\phi_k,\ \phi\in\{0,1\}^{d_m}\right\}\ ,$$
if the positive number $r$ is small enough so that $\mathcal{C}_m(\theta',r)\subset \Theta^+$. For any $\theta'\in\Theta_m^{+,\text{\emph{iso}}}$, we analogously define the hypercubes $\mathcal{C}^{\text{\emph{iso}}}_m\left(\theta',r\right)$ using a basis $\left(\Psi^{\text{\emph{iso}}}_{i_1,j_1},\ldots, \Psi^{\text{\emph{iso}}}_{i_{d_m},j_{d_m}}\right)$.
\end{defi}

\begin{prte}\label{variance_minimax_global}
Let $m$ be a  model in $\mathcal{M}_1\setminus\emptyset$ whose dimension $d_m$ is smaller than $p\sqrt{n}$. Then, for any estimator $\widehat{\theta}$, 
\begin{eqnarray}\label{minoration_parametrique}
\sup_{\theta\in \Theta_{m}^+}\mathbb{E}_{\theta}\left[l\left(\widehat{\theta},\theta\right)\right]\geq \sup_{\theta\in \Theta_{m,2}^+}\mathbb{E}_{\theta}\left[l\left(\widehat{\theta},\theta\right)\right]\geq L\sigma^2\frac{d_m}{np^2}\ .
\end{eqnarray}
Let $\theta'$ be an element of $\Theta_m^+$ that satisfies $(\mathbb{H}_2)$. For any estimator $\widehat{\theta}$ of $\theta$, 
\begin{eqnarray}\label{minimax_local}
 \sup_{\theta\in \text{\emph{Co}}\left[\mathcal{C}_m\left(\theta',(1-\|\theta'\|_1)/\sqrt{np^2}\right)\right]}\mathbb{E}_{\theta}\left[l\left(\widehat{\theta},\theta\right)\right]\geq L\sigma^2\varphi^2_{\text{\emph{min}}}\left[I_{p^2}-C(\theta')\right]\frac{d_m}{np^2}\ ,
\end{eqnarray}
where $\text{\emph{Co}}\left[\mathcal{C}_m\left(\theta',r\right)\right]$ denotes the convex hull of $\mathcal{C}_m\left(\theta',r\right)$.
\end{prte}
An analogous result holds for isotropic hypercubes.
The first bound (\ref{minoration_parametrique}) means that for any estimator $\widehat{\theta}$, the supremum of the risks $\mathbb{E}_{\theta}[l(\widehat{\theta}_{m,\rho_1},\theta)]$ over $\Theta_m^+$ is larger than $\sigma^2d_m/(np^2)$ (up to some numerical constant). This rate $\sigma^2d_m/(np^2)$ is achieved by the CLS estimator by Theorem \ref{mainthrm}.

The second lower bound (\ref{minimax_local}) is of independent interest. It implies that in a small neighborhood of $\theta'$ the risk $\mathbb{E}_{\theta}[l(\widehat{\theta}_{m,\rho_1},\theta)]$ is larger than $\sigma^2\varphi^2_{\text{min}}[I_{p^2}-C(\theta')] d_m/(np^2)$. This confirms the lower bound (\ref{minorationas}) of Corollary \ref{minoration_variance_asymptotique} in a nonasymptotic way. Indeed, these two expressions match up to a factor $\varphi_{\text{min}}[I_{p^2}-C(\theta')] $. This difference comes from the fact that the lower bound (\ref{minimax_local}) holds for any estimator $\widehat{\theta}$. Bound (\ref{minimax_local}) is sharp in the sense that the maximum likelihood estimator $\widehat{\theta}_{m_1}^{\text{iso},mle}$ of isotropic GMRF in $m_1$ exhibits an asymptotic risk of order $\sigma^2\varphi^2_{\text{min}}[I_{p^2}-C(\theta)]/(np^2)$ for the parameter $\theta$ studied in Example \ref{variance_limite_alpha}. It is shown using the methodology introduced in the proof of Example \ref{variance_limite_alpha}.
We now state that $\widetilde{\theta}_{\rho}$ is adaptive to the sparsity of  $m$.
\begin{cor}\label{adaptation_sparsite}
Considering $K\geq K_0$, $\rho_1\geq 2$, $\rho_2>2$ and a collection $\mathcal{M}\subset\mathcal{M}_1$, we define the estimator $\widetilde{\theta}_{\rho_1}$ with the penalty $\pen(m)=K\sigma^2\rho_1^2\rho_2\frac{d_m}{np^2}$. For any non empty model $m$,   
\begin{eqnarray}\label{adaptation_sparsite_equation}
 \sup_{\theta\in\Theta^+_{m,\rho_1}\cap\mathcal{U}(\rho_2)}\mathbb{E}_{\theta}\left[l\left(\widetilde{\theta}_{\rho_1},\theta\right)\right]\leq L\left(K,\rho_1,\rho_2\right)\inf_{\widehat{\theta}}\sup_{\theta\in\Theta^+_{m,\rho_1}\cap\mathcal{U}(\rho_2)}\mathbb{E}\left[l\left(\widehat{\theta},\theta\right)\right]\ ,
\end{eqnarray}
where $\mathcal{U}(\rho_2)$ is defined in (\ref{definition_U}).
\end{cor}
A similar result holds for $\widetilde{\theta}_{\rho_1}^{\text{iso}}$ and $\Theta^{+,\text{iso}}_{m,\rho_1}$. Corollary \ref{adaptation_sparsite} is nonasymptotic and applies for any $n$ and any $p$. 
If $\theta$ belongs to some model $m$, then the optimal risk from a minimax point of view is of order $\frac{d_m}{np^2}$. In practice, we do not know the true model $m$. Nevertheless, the procedure simultaneously achieves the minimax rates for all supports $m$ possible.
This means that $\widetilde{\theta}_{\rho_1}$ reaches this minimax rate $\frac{d_m}{np^2}$ \emph{without} knowing in advance the true model $m$. 

The procedure is not adaptive to the smallest and the largest eigenvalue of $(I_{p^2}-C(\theta))$ which correspond to $\rho_1$ and $\rho_2$. Indeed, the constant $L\left(K,\rho_1,\rho_2\right)$ depends on $\rho_1$ and $\rho_2$. We are not aware of any other covariance estimation procedure which is really adaptive the smallest and the largest eigenvalue of the matrix.~\\~\\
Finally, $\widetilde{\theta}_{\rho_1}$ exhibits the same adaptive properties with respect to the Frobenius norm.
\begin{cor}\label{adaptation_sparsite_Frobenius}
Under the same assumptions as Corollary \ref{adaptation_sparsite},   
\begin{eqnarray*}
 \sup_{\theta\in\Theta^+_{m,\rho_1}\cap\mathcal{U}(\rho_2)}\mathbb{E}_{\theta}\left[\|C(\widetilde{\theta}_{\rho_1})-C(\theta)\|_F^2\right]\leq L\left(K,\rho_1,\rho_2\right)\inf_{\widehat{\theta}}\sup_{\theta\in\Theta^+_{m,\rho_1}\cap\mathcal{U}(\rho_2)}\mathbb{E}\left[\|C(\widehat{\theta})-C(\theta)\|_F^2\right]\ .
\end{eqnarray*}
\end{cor}

\begin{proof}[Proof of Corollary \ref{adaptation_sparsite_Frobenius}]
As in the proof of Corollary \ref{corollaire_Frobenius}, we observe that 
$$\|C(\theta_1)-C(\theta_2)\|_F\geq \frac{p^2\rho_1}{\sigma^2} l(\theta_1,\theta_2)\ ,$$ if $\theta$ satisfies Assumption $(\mathbb{H}_1)$.
We conclude by applying Proposition \ref{variance_minimax_global} and Corollary \ref{corollaire_Frobenius}.
\end{proof}

\subsection{Adapting to the decay of the bias}\label{section_adaptive_ellipsoids}

In this section, we prove that the estimator $\widetilde{\theta}_{\rho_1}$ is adaptive to a range of sets that we call \emph{pseudo-ellipsoids}.
\begin{defi}[Pseudo-ellipsoids]\label{definition_ellipsoide}
Let $(a_j)_{1\leq j \leq \text{Card}(\mathcal{M}_1)}$ be a non-increasing sequence of
positive numbers. Then,  $\theta\in\Theta^+$ belongs to the \emph{pseudo-ellipsoid}	
$\mathcal{E}(a)$ if and only if
\begin{eqnarray}
  \sum_{i=1}^{\text{Card}(\mathcal{M}_1)}\frac{\emph{var}_{\theta}\left(X{\scriptstyle[0,0]}|X_{\mathcal{N}(m_{i-1})}\right)-\emph{var}_{\theta}\left(X{\scriptstyle[0,0]}|X_{\mathcal{N}({m_{i}})}\right)}{a_i^2}& \leq &
  1\ .\label{def_ellipsoide_2}
\end{eqnarray}
\end{defi}
Condition (\ref{def_ellipsoide_2}) measures how fast $\var_{\theta}(X{\scriptstyle [0,0]}|X_{\mathcal{N}(m_{i})})$ tends to $\var_{\theta}(X{\scriptstyle[0,0]}|X_{\Lambda\setminus\{(0,0)\}})$.
Suppose that Assumption $(\mathbb{H}_2)$  defined in (\ref{def_h2}) is fulfilled. By Corollary \ref{definition_thetam}, $\var_{\theta}\left(X{\scriptstyle[0,0]}|X_{\mathcal{N}(m_{i})}\right)$ is  the sum of  $l(\theta_{m_i},\theta)$ and $\sigma^2$ and  Condition (\ref{def_ellipsoide_2}) is equivalent to 
\begin{eqnarray}\label{def_ellipsoide_3}
\sum_{i=1}^{\text{Card}(\mathcal{M}_1)}\frac{l\left(\theta_{m_{i-1}},\theta)-l(\theta_{m_{i}},\theta\right)}{a_i^2}\leq 1\ .
\end{eqnarray}
Hence, the sequence ($a_i$) gives some condition on the \emph{rate of decay} of the bias when the dimension of the model increases. These sets $\mathcal{E}(a)$ are not true ellipsoids.  Nevertheless, one may consider them as counterparts of the classical ellipsoids studied in the fixed design Gaussian regression framework (see for instance \cite{massartflour} Sect.4.3). 

To prove adaptivity, we shall need the equivalence between Conditions   (\ref{def_ellipsoide_2}) and (\ref{def_ellipsoide_3}). This equivalence holds if   $\var_{\theta}\left(X{\scriptstyle[0,0]}|X_{\mathcal{N}(m_{i})}\right)$ decomposes as $l(\theta_{m_i},\theta)+\sigma^2$, for any model $m\in\mathcal{M}_1$. As mentioned earlier, Assumption $(\mathbb{H}_2)$ is sufficient (but not necessary) for this property to hold. This is why we restrict ourselves to study sets of the type $\mathcal{E}(a)\cap \mathcal{B}_1(0_p,1)$.
We shall also perform the following assumption on the ellipsoids $\mathcal{E}(a)$ 
\begin{eqnarray*}
(\mathbb{H}_a): \hspace{4cm}  a_i^2\leq \frac{\sigma^2}{d_{m_i}},\text{ for any $1\leq i\leq \left|\mathcal{M}_1\right|$}\ .\hspace{5cm}
\end{eqnarray*}
It essentially means that the sequence $(a_i)$ converges fast enough towards $0$. For instance, all the sequences $a_i=\sigma(d_{m_i})^{-s}$ with $s\geq 1/2$ satisfy $(\mathbb{H}_a)$.

\begin{prte}\label{minoration_ellipsoid} 
Under Assumption $(\mathbb{H}_a)$, the minimax rate of estimation on $\mathcal{E}(a)\cap\mathcal{B}_1(0_p,1)\cap\mathcal{U}(2)$ is lower bounded by
\begin{eqnarray}\label{minoration_ellipsoid_equation}
\inf_{\widehat{\theta}}\sup_{\theta \in \mathcal{E}(a)\cap\mathcal{B}_1(0_p,1)\cap\mathcal{U}(2)}\mathbb{E}_{\theta}\left[l\left(\widehat{\theta},\theta\right)\right] \geq L
\sup_{1\leq i\leq \text{Card}(\mathcal{M}_1)}\left(a_i^2\wedge \sigma^2\frac{d_{m_i}}{np^2}\right)\ .
\end{eqnarray}
\end{prte}

This lower bound is analogous to the minimax rate of estimation for ellipsoids in the Gaussian sequence model.
Gathering Theorem \ref{mainthrm} and Proposition \ref{minoration_ellipsoid} enables to derive adaptive properties for $\widetilde{\theta}_{\rho_1}$.

\begin{prte}\label{prte_adaptatif}
Considering $K\geq K_0$, $\rho_1\geq2$, $\rho_2>2$ and the collection $\mathcal{M}_1$, we define the estimator $\widetilde{\theta}_{\rho_1}$ with the penalty $\pen(m)=K\sigma^2\rho_1^2\rho_2\frac{d_m}{np^2}$. 
For any ellipsoid $\mathcal{E}(a)$ that satisfies $(\mathbb{H}_a)$ and such that $a_1^2 \geq 1/(np^2)$, the estimator $\widetilde{\theta}_{\rho_1}$ is minimax  over the set $\mathcal{E}(a)\cap\mathcal{B}_1(0_p,1)\cap\mathcal{U}(\rho_2)$:
\begin{eqnarray}\label{eq_adaptatif}
\sup_{\theta \in \mathcal{E}(a)\cap\mathcal{B}_1(0_p,1)\cap\mathcal{U}(\rho_2)}\mathbb{E}_{\theta}\left[l\left(\widetilde{\theta}_{\rho_1},\theta\right) \right]\leq
L(K,\rho_1,\rho_2)\inf_{\widehat{\theta}} \sup_{\theta\in
  \mathcal{E}(a)\cap\mathcal{B}_1(0_p,1)\cap\mathcal{U}(\rho_2)} \mathbb{E}_{\theta}\left[l\left(\widehat{\theta},\theta\right)\right] \ .
\end{eqnarray}
\end{prte}

Let us first illustrate this result. We have mentioned earlier, that Assumption $(\mathbb{H}_a)$ is satisfied for all sequences $a_i=\sigma(d_{m_i})^{-s}$ with $s\geq 1/2$. We note $\mathcal{E}'(s)$ such a pseudo-ellipsoid. By Propositions \ref{minoration_ellipsoid} and \ref{prte_adaptatif}, the minimax rate over \emph{one} pseudo ellipsoid $\mathcal{E}'(s)$ is $\sigma^2(np^2)^{-2s/(1+2s)}$. The larger $s$ is, the faster the minimax rates is. The estimator $\widetilde{\theta}_{\rho_1}$ achieves simultaneously the rate $\sigma^2(np^2)^{-2s/(1+2s)}$ for all $s\geq 1/2$. Consequently, $\widetilde{\theta}_{\rho_1}$ is adaptive to the rate $s$ of decay of the bias: it achieves the optimal rates without knowing $s$ in advance.

Let us further comment Proposition  \ref{prte_adaptatif}. By (\ref{eq_adaptatif}), the estimator $\widetilde{\theta}_{\rho_1}$ is  adaptive over $\mathcal{E}(a)\cap\mathcal{B}_1(0_p,1)\cap\mathcal{U}(\rho_2)$ for all sequences $(a)$ such that  $(\mathbb{H}_a)$ is satisfied and such that $a_1^2\geq 1/(np^2)$. Again, the result applies for any $n$ and any $p$.
The condition $a_1^2\geq 1/(np^2)$ is classical. It ensures that the pseudo-ellipsoid $\mathcal{E}(a)$ is not degenerate, i.e. that the minimax rates of estimation is not smaller than $\sigma^2/(np^2)$. We have explained earlier that we restricts ourselves to  parameters $\theta$ in $\mathcal{B}_1(0_p,1)$ only because this enforces the equivalence between (\ref{def_ellipsoide_2}) and (\ref{def_ellipsoide_3}). In contrast, the hypothesis $\varphi_{\text{max}}(\Sigma)\leq \sigma^2\rho_2$ is really necessary because we fail to be adaptive to $\rho_2$.

\begin{cor}\label{adaptivity_frobenius}
Under Assumption $(\mathbb{H}_a)$, the minimax rate of estimation over $\mathcal{E}(a)\cap\mathcal{U}(2)\cap\mathcal{B}_1(0_p,1)$ is lower bounded by
\begin{eqnarray*}
\inf_{\widehat{\theta}}\sup_{\theta \in \mathcal{E}(a)\cap\mathcal{B}_1(0_p,1)\cap\mathcal{U}(2)}\mathbb{E}_{\theta}\left[\|C(\widehat{\theta})-C(\theta)\|_F^2\right] \geq L
\sup_{1\leq i\leq \text{Card}(\mathcal{M}_1)}\left(a_i^2p^2\wedge \frac{d_{m_i}}{n}\right)\ .
\end{eqnarray*}
Under the same assumptions as Proposition \ref{prte_adaptatif}, 
\begin{eqnarray*}
\sup_{\theta \in \mathcal{E}(a)\cap\mathcal{B}_1(0_p,1)\cap\mathcal{U}(\rho_2)}\mathbb{E}_{\theta}\left[\|C(\widehat{\theta})-C(\theta)\|_F^2 \right]\leq
L(K,\rho_1,\rho_2)\inf_{\widehat{\theta}} \sup_{\theta\in
  \mathcal{E}(a)\cap\mathcal{B}_1(0_p,1)\cap\mathcal{U}(\rho_2)} \mathbb{E}_{\theta}\left[\|C(\widehat{\theta})-C(\theta)\|_F^2\right] \ .
\end{eqnarray*}
\end{cor}

\begin{proof}[Proof of Corollary \ref{adaptivity_frobenius}]
As in the proof of Corollary \ref{corollaire_Frobenius}, we observe that $$\|C(\theta_1)-C(\theta_2)\|_F\geq p^2[\varphi_{\text{max}}(\Sigma)]^{-1} l(\theta_1,\theta_2)\geq \frac{p^2}{\rho_2\sigma^2} l(\theta_1,\theta_2)\ ,$$
$$\|C(\theta_1)-C(\theta_2)\|_F\leq p^2[\varphi_{\text{min}}(\Sigma)]^{-1}l(\theta_1,\theta_2)\leq p^2 \frac{\varphi_{\text{max}}[I_{p^2}-C(\theta)]}{\sigma^2}l(\theta_1,\theta_2)\leq \frac{\rho_2p^2}{\sigma^2}l(\theta_1,\theta_2)\ ,$$
 if $\theta\in\mathcal{B}_1(0_p,1)\cap\mathcal{B}_{\text{op}}(\rho_2)$.
We conclude by applying Proposition \ref{minoration_ellipsoid} and Proposition \ref{prte_adaptatif}.
\end{proof}

Again, $\widetilde{\theta}_{\rho_1}$ satisfies the same minimax properties with respect to the Frobenius norm. All these properties easily extend to isotropic fields if one defines the corresponding sets $\mathcal{E}^{\text{iso}}(a)\cap \mathcal{B}_1(0_p,1)\cap\mathcal{U}(\rho_2)$ of isotropic GMRFs.

\section{Discussion}\label{section_discussion}

\subsection{Comparison with maximum likelihood estimation}\label{section_mle}

%%% Reecrire
Let us first compare the computational cost the CLS estimation method and the maximum likelihood estimator (MLE).
 For toroidal lattices, fast algorithms  based on two-dimensional fast-Fourier transformation (see for instance  \cite{rue02}) allow to compute the MLE as fast as the CLS estimator. More details on the computation of the CLS estimators for toroidal lattices are given in \cite{verzelen_gmrf_appli} Sect.2.3. 
When the lattice is not a torus, the MLE becomes intractable because it involves the optimization of a determinant of size $p^2$. In contrast, the CLS criterion $\gamma_{n,p}(.)$ defined in (\ref{definition_gamman}) is a quadratic function of $\theta$. Consequently, CLS estimators are still computationally amenable. We extend our model selection to non-toroidal lattices in \cite{verzelen_gmrf_appli}.\\

Let us compare the risk of CLS estimators and MLE.
Given a small dimensional model $m$, the risk of the \emph{parametric} CLS estimator and the \emph{parametric} MLE have been compared from an asymptotic point of view (\cite{guyon95} Sect.4.3). It is generally accepted  (see for instance Cressie \cite{cressie} Sect. 7.3.1) and that \emph{parametric} CLS estimators are almost as efficient as parametric MLE for the major part of the parameter spaces $\Theta_{m}^+$. We have non-asymptotically assessed this statement in Proposition  \ref{variance_minimax_global} by minimax arguments. Nevertheless, for some parameters $\theta$ that are close to the border of $\Theta_{m}^+$, Kashyap and Chellappa \cite{kashyap83} have pointed out that CLS estimators are  less efficient than MLE. If we have proved nonasymptotic bounds for CLS-based model selection method, we are not aware of any such result for model selection procedures based on MLE.

\subsection{Concluding remarks}

We have developed a model selection procedure for choosing the neighborhood of a GMRF. In Theorem \ref{mainthrm}, we have proven  a nonasymptotic upper bound for the risk of the estimator $\widetilde{\theta}_{\rho_1}$ with respect to the prediction error $l(.,.)$. Under Assumption $(\mathbb{H}_1)$, this bound is shown to be optimal from an asymptotic point of view if the support of $\theta$ belongs to one of the models in the collection. If Assumption $(\mathbb{H}_2)$  is fulfilled, we are able to obtain an oracle type inequality for $\widetilde{\theta}_{\rho_1}$. Moreover, $\widetilde{\theta}_{\rho_1}$ is minimax adaptive to the sparsity of $\theta$ under $(\mathbb{H}_1)$. Finally, it simultaneously achieves the minimax rates of estimation over a large class of sets $\mathcal{E}(a)$ if $(\mathbb{H}_2)$ holds. Some of these properties still hold if we use the Frobenius loss function. The case of isotropic Gaussian fields is handled similarly.\\

However, in the oracle inequality (\ref{inegalite_oracle_asymptotique}) and in the minimax bounds (\ref{adaptation_sparsite_equation}) and (\ref{eq_adaptatif}), we either perform an assumption on the $l_1$ norm of $\theta$ or on the smallest eigenvalue of $(I_{p^2}-C(\theta))$. When $\|\theta\|_1$ tends  to one or $\varphi_{\text{min}}[I_{p^2}-C(\theta)]$ tends  to $0$, there is a distortion between the upper bound $\mathbb{E}_{\theta}[l(\widetilde{\theta}_{\rho_1},\theta)]$ provided by Theorem \ref{mainthrm} and the lower bounds given by Corollary \ref{minoration_variance_asymptotique} or Proposition \ref{variance_minimax_global}. This limitation seems intrinsic to our penalization method which is linear with  respect to the dimension, whereas the asymptotic variance term $\mathbb{E}_{\theta}[l(\widehat{\theta}_{m,\rho_1},\theta)]$ depends in a complex way on the dimension of the model $m$ and on the target $\theta$. In our opinion, achieving adaptivity with respect to the smallest eigenvalue of $(I_{p^2}-C(\theta))$ (or equivalently the largest value of $\Sigma$) would require a different penalization technique. Nevertheless, we are not aware of any procedure in a covariance estimation setting that is adaptive to the largest eigenvalues of $\Sigma$.\\

So far, we have provided an estimation procedure for $(I_{p^2}-C(\theta))=\sigma^2\Sigma^{-1}$. If we aim at estimating the precision matrix $\Sigma^{-1}$, we also have to take into account the quantity $\sigma^2$. It is natural to estimate it by $\widetilde{\sigma}^2:=\gamma_{n,p^2}(\widetilde{\theta}_{\rho_1})$ as done for instance by Guyon in \cite{guyon95} Sect.4.3 in the parametric setting. Then, we obtain the estimate $\widetilde{\Sigma^{-1}}:=\widetilde{\sigma}^2(I_{p^2}-C(\widetilde{\theta}_{\rho_1})$. It is of interest to study the adaptive properties of this estimator with respect to loss functions such as the Frobenius or operator norm as done in \cite{rothman07} in the non-stationary setting. Nevertheless, let us mention that the matrix $\widetilde{\Sigma^{-1}}$ is not necessarily invertible since the estimator $\widetilde{\theta}_{\rho_1}$ belongs to the closure of $\Theta^+$.\\

The choice of the quantity $\rho_1$ is  problematic.  On the one hand, $\rho_1$ should be large enough so that Assumption $(\mathbb{H}_1)$ is fulfilled. On the other hand, a large value of $\rho_1$ yields worse bounds in Theorem \ref{mainthrm}.
Moreover, the largest eigenvalue of $(I_{p^2}-C(\theta))$ is unknown in practice, which makes more difficult the choice of $\rho_1$.
We see two possible answers to this issue:
\begin{itemize}
 \item First, moderate values of $\rho_1$ are sufficient to enforce $(\mathbb{H}_1)$ if the target $\theta$ is sparse as illustrated in Table \ref{tableau_comparaison_rho}.
\item Second,  we believe that the bounds for the risk are pessimistic with respect to $\rho_1$. A future direction of research is to derive risk bounds for $\widetilde{\theta}_{\rho_1}$ with $\rho_1=+\infty$. In \cite{verzelen_gmrf_appli}, we illustrate that such a procedure gives rather good results in practice.
\end{itemize}

In Theorem \ref{mainthrm}, we only provide a lower bound of the penalty so that the procedure performs well. However, this bound depends on the largest eigenvalue of $\Sigma$ which is seldom known in practice and we did not give any advice for 
choosing a ``reasonable'' constant $K$ in practice. This is why we introduce in \cite{verzelen_gmrf_appli} a data-driven method based on the \emph{slope heuristics} of Birg\'e and Massart \cite{massart_pente}  for calibrating  the penalty. We also provide numerical evidence of its performances on simulated data. For instance, the procedure outperforms variogram-based methods for estimating Mat\'ern correlations.\\

We have mentioned in the introduction that the toroidal assumption for the lattice is somewhat artificial in several applications. Nevertheless, we needed to neglect the edge effects in order to derive non asymptotic properties for $\widetilde{\theta}_{\rho_1}$ as in Theorem \ref{mainthrm}. In practice, it is often more realistic to suppose that we observe a small window of a Gaussian field defined on the whole plane $\mathbb{Z}^2$. The previous nonasymptotic properties do not extend to this new setting. Nevertheless, Lakshman and Derin have shown in \citep{Lakshman93} that there is no phase transition within the valid parameter space for GMRFs defined on the plane $\mathbb{Z}^2$. In short, this implies that the distribution of a field observed in a fixed window of a GMRF does not asymptotically depend on the bound condition. Therefore, it is reasonable to think that our estimation procedure performs well if it was adapted  to this new setting. In \cite{verzelen_gmrf_appli}, we describe such an extension and we provide numerical evidence of its performances.

\subsection{Possible extensions}

In many statistical applications stationary Gaussian fields (or Gaussian Markov random fields) are not directly observed. For instance, Aykroyd \cite{aykroyd} or Dass and Nair \cite{dass03} use compound Gaussian Markov random fields to account for non stationarity and steep variations. The wavelet transform has emerged as a powerful tool in image analysis. the wavelet coefficients of an image are sometimes modeled using hidden Markov models~\cite{crouse98,portilla}. More generally, the success of the GMRF is mainly due to the use of hierarchical models involving latent GMRFs~\cite{inla}. The study and the implementation of our penalization strategy for selecting the complexity of the latent Markov models is an interesting direction of research.

\section{Proofs}\label{section_Preuves}

\subsection{A concentration inequality}\label{preuves_concentration}

In this section, we prove a new concentration inequality for suprema of Gaussian chaos of order 2. It will be useful for proving Theorem \ref{mainthrm}.

\begin{prte}\label{lemmechaos}
Let $F$ be a compact set of symmetric matrices of size $r$, $(Y^1,\ldots,Y^n)$ be
a $n$-sample of a standard Gaussian vector of size $r$, and $Z$ be the random variable defined by
$$Z:= \sup_{R \in F}tr\left[R(\overline{YY^*}-I_r)\right]\ .$$
Then 
\begin{eqnarray} \label{inegalite_concentration}
\mathbb{P}(Z \geq \mathbb{E}(Z) + t) \leq \exp\left[-\left( \frac{t^2}{L_1\mathbb{E}(W)}\bigwedge \frac{t}{L_2B}\right)\right],
\end{eqnarray}
where the quantities $B$ and $W$ are such that
\begin{eqnarray} 
B & := & \frac{2}{n} \sup_{R \in F}\varphi_{\text{max}}(R) \nonumber\\
W & := & \frac{4}{n} \sup_{R \in F}tr(R\overline{YY^*}R')\ .\nonumber
\end{eqnarray}
\end{prte}

The main argument of this proof is to transfer a deviation inequality 
 for suprema of Rademacher chaos of order 2 to suprema of Gaussian Chaos. Talagrand \cite{Talagrand96} has first given in Theorem 1.2 a concentration inequality for such suprema of Rademacher chaos. Boucheron \emph{et al.} \cite{bousquet05} have recovered the upper bound applying a new methodology based on the entropy method. We  adapt their proof to consider non-necessarily homogeneous chaos of order $2$. More details are found in the technical appendix \cite{technical}.

\subsection{Proof of Theorem \ref{mainthrm}}\label{section_preuve_mainthrm}

\begin{proof}[Proof of Theorem \ref{mainthrm}]
We only consider the case of anisotropic estimators. The proofs and lemma are analogous for isotropic estimators.
 We first fix a model $m\in\mathcal{M}$.
By definition, the model $\widehat{m}$ satisfies
$$\gamma_{n,p}(\widetilde{\theta}_{\rho_1})+\pen(\widehat{m}) \leq \gamma_{n,p}(\theta_{m,\rho_1})+ \pen(m)\ .$$
For any $\theta'\in\Theta^+$, $\overline{\gamma}_{n,p}(\theta')$ stands for the difference between
$\gamma_{n,p}(\theta')$ and its expectation $\gamma(\theta')$. Then,
the previous inequality turns into
\begin{eqnarray*}
\gamma(\widetilde{\theta}_{\rho_1}) \leq  \gamma(\theta_{m,\rho_1}) +\overline{\gamma}_{n,p}(\theta_{m,\rho_1}) -\overline{\gamma}_{n,p}(\widetilde{\theta}_{\rho_1}) +\pen(m) - \pen(\widehat{m})\ .
\end{eqnarray*}
Subtracting the quantity $\gamma(\theta)$ to both sides of this inequality yields 
\begin{eqnarray}\label{preuve2}
l(\widetilde{\theta}_{\rho_1},\theta) \leq l(\theta_{m,\rho_1},\theta)+ \overline{\gamma}_{n,p}(\theta_{m,\rho_1}) -\overline{\gamma}_{n,p}(\widetilde{\theta}_{\rho_1})  + \pen(m)-\pen(\widehat{m})\ .
\end{eqnarray}
The proof is based on the control of the random variable $\overline{\gamma}_{n,p}(\theta_{m,\rho_1}) -\overline{\gamma}_{n,p}(\widetilde{\theta}_{\rho_1})$.

\begin{lemma}\label{thrm_argument_principal}
 For any positive number $\alpha$, $\xi$, and $\delta>1$ the event  $\Omega_{\xi}$ defined by 
\begin{eqnarray*}
 \Omega_{\xi} = \left\{\begin{array}{r}

\overline{\gamma}_{n,p}\left(\theta_{m,\rho_1}\right)-\overline{\gamma}_{n,p}\left(\widetilde{\theta}_{\rho_1}\right)\leq \frac{1}{\sqrt{\delta}}l\left(\widetilde{\theta}_{\rho_1},\theta\right)+\frac{\sqrt{\delta}}{\sqrt{\delta}-1}l\left(\theta_{m,\rho_1},\theta\right)\\ +\frac{K_0\delta^2\rho_1^2\varphi_{\text{max}}(\Sigma)}{np^2}\left[(1+\alpha/2)\left(d_{m}+d_{\widehat{m}}\right)+\frac{\xi^2}{\delta-1}\right]
\end{array}
\right\} \ , 
\end{eqnarray*}
 satisfies 
\begin{eqnarray*}
 \mathbb{P}(\Omega_{\xi}^c)\leq\exp\left\{-L_1\xi\left[\frac{\alpha}{\sqrt{1+\alpha/2}}\wedge \sqrt{n}\right] \right\} \sum_{m'\in\mathcal{M}}\exp\left\{-
 L_2\sqrt{d_{m'}}\left(\frac{\alpha}{\sqrt{1+\alpha/2}}\wedge
 \frac{\alpha^2}{1+\alpha/2}\right)\right\}\ .
\end{eqnarray*}
\end{lemma}
A similar lemma holds in the isotropic case.
In particular, we choose $\alpha=(K-K_0)/K_0$ and $\delta =\sqrt{(1+\alpha)/(1+\alpha/2)}$. Lemma \ref{thrm_argument_principal} implies that on the event $\Omega_\xi$, 
\begin{eqnarray*}
 \overline{\gamma}_{n,p}\left(\theta_{m,\rho_1}\right)-\overline{\gamma}_{n,p}\left(\widetilde{\theta}_{\rho_1}\right)& \leq & \frac{1}{\sqrt{\delta(\alpha)}}l\left(\widetilde{\theta}_{\rho_1},\theta\right)+\frac{\sqrt{\delta(\alpha)}}{\sqrt{\delta(\alpha)}-1}l\left(\theta_{m,\rho_1},\theta\right)+ \pen(m)\\ &+ &\pen(\widehat{m}) +\frac{K_0\xi^2\delta(\alpha)^2\rho_1^2\varphi_{\text{max}}(\Sigma)}{np^2\left(\delta(\alpha)-1\right)}\ .
\end{eqnarray*}
Thus, gathering this bound with inequality (\ref{preuve2}) yields
\begin{eqnarray*} 	
\frac{\delta(\alpha)^{1/2}-1}{\delta(\alpha)^{1/2}}l(\widetilde{\theta}_{\rho_1},\theta) &\leq& \left[1+\delta(\alpha)^{-1/2}(\delta(\alpha)^{1/2}-1)^{-1}\right]l(\theta_{m,\rho_1},\theta) + 2 \pen(m)\\& +&
 \frac{K_0\xi^2\rho_1^2\varphi_{\text{max}}(\Sigma) \delta(\alpha)^2}{np^2(\delta(\alpha)-1)}\ , \nonumber
\end{eqnarray*}
with probability larger than $1-\mathbb{P}(\Omega_{\xi})$.
Integrating this inequality with respect to $\xi>0$ leads to
\begin{eqnarray}\label{equiintegrable}
\frac{\delta(\alpha)^{1/2}-1}{\delta(\alpha)^{1/2}}\mathbb{E}_{\theta}\left[l\left(\widetilde{\theta}_{\rho_1},\theta\right)\right]
\leq \nonumber \left[1+\delta(\alpha)^{-1/2}\left(\delta(\alpha)^{1/2}-1\right)^{-1}\right]l\left(\theta_{m,\rho_1},\theta\right) + \\2 \pen(m)+
\frac{ \delta(\alpha)^2L(\alpha)}{\left(\delta(\alpha)-1\right)\left[ \frac{\alpha^2}{1+\alpha/2}\wedge n\right]}\frac{\rho_1^2\varphi_{\text{max}}(\Sigma)}{np^2}\ .
\end{eqnarray}
We upper bound $[(\alpha^2/(1+\alpha/2)) \wedge n]^{-1}$ by.$[(\alpha^2/(1+\alpha/2)) \wedge 1]^{-1}$.
Since $\alpha= \frac{K-K_0}{K_0}$, it follows that 
\begin{eqnarray*}
\mathbb{E}_{\theta}\left[l\left(\widetilde{\theta}_{\rho_1},\theta\right)\right] \leq L_1(K)\left[l\left(\theta_{m,\rho_1},\theta\right)+\pen(m) \right] + L_2\left(K\right)\frac{\rho_1^2\varphi_{\text{max}}(\Sigma)}{np^2}\ ,
\end{eqnarray*}
 Taking the infimum over the models $m\in \mathcal{M}$ allows to conclude.
\end{proof}\vspace{0.5cm}

\begin{proof}[Proof of Lemma \ref{thrm_argument_principal}]
Throughout this proof, it is more convenient to express the
quantities $\overline{\gamma}_{n,p}(.)$ and $l(.)$ in terms of covariance and precision matrices. Thanks to Equation  (\ref{definition_gamma}), we also provide a matricial expression for $\gamma(.)$ :
\begin{eqnarray}\label{definition_gamma_alternative}
\gamma(\theta') = \frac{1}{p^2}tr\left[\left(I-C(\theta')\right)\Sigma\left(I-C(\theta')\right)\right] \ .
\end{eqnarray}
Gathering identities  (\ref{definition_gamma_alternative}) and (\ref{definition_gamman_alternative}), we get
\begin{eqnarray*}
\overline{\gamma}_{n,p}(\theta_{m,\rho_1}) -\overline{\gamma}_{n,p}(\widetilde{\theta}_{\rho_1}) & = &
\frac{1}{p^2}tr\left[\left(\left[I_{p^2}-C(\theta_{m,\rho_1})\right]^2
  -\left[I_{p^2}-C(\widetilde{\theta}_{\rho_1})\right]^2\right)\left(\overline{\bf X^vX^{v*}}-\Sigma\right)\right]\ . \label{preuve_matricielle} 
\end{eqnarray*}

Since the matrices $\Sigma$, $(I_{p^2}-C(\theta_{m,\rho_1}))$, and $(I_{p^2}-C(\widetilde{\theta}_{\rho_1}))$ correspond to 
covariance or precision matrices of stationary fields on the two dimensional torus, they
are symmetric block circulant. By Lemma \ref{codiagonalisation}, they are jointly diagonalizable
in the same orthogonal basis. In the sequel, $P$ stands for an orthogonal matrix associated
to this basis. Then, the matrices $C(\theta_{m,\rho_1})$,
$C(\widetilde{\theta}_{\rho_1})$, and $\Sigma$ respectively decompose in
$$C(\theta_{m,\rho_1})=P^*D({\theta_{m,\rho_1}})P,\text{     }C(\widetilde{\theta}_{\rho_1})=P^*D({\widetilde{\theta}_{\rho_1}})P, \text{      }\Sigma=P^*D_{\Sigma}P,$$
where the matrices $D(\theta_{m,\rho_1})$, $D(\widetilde{\theta}_{\rho_1})$, and $D_{\Sigma}$ are
diagonal. Let the $p^2\times n$ matrix ${\bf Y}$ be defined by
${\bf Y}:=\sqrt{\Sigma^{-1}}{\bf X^v}$. Clearly, the components of ${\bf Y}$ follow independent standard normal distributions. Gathering these new notations, we get 
\begin{eqnarray}
\lefteqn{\overline{\gamma}_{n,p}(\theta_{m,\rho_1})-\overline{\gamma}_{n,p}(\widetilde{\theta}_{\rho_1})
\, =} \nonumber &  &\\ & & \frac{1}{p^2}tr\left[\left(\left[I_{p^2}-D({\theta_{m,\rho_1}})\right]^2-\left[I_{p^2}-D({\widetilde{\theta}_{\rho_1}})\right]^2\right)D_{\Sigma}\left(\overline{\bf YY^*}-I_{p^2}\right)\right]\ . \label{borne}
\end{eqnarray}
Except $\overline{\bf YY^*}$ all the matrices in this last expression are diagonal and we may therefore 
commute them in the trace.\vspace{0.5cm}~\\
Let $<.,.>_{\mathcal{H}}$ and $<.,.>_{\mathcal{H}'}$ be two inner products in the space of
square matrices of size $p^2$ respectively defined by 
$$<A,B>_{\mathcal{H}}:=\frac{tr(A^*\Sigma B)}{p^2}\text{ and } <A,B>_{\mathcal{H}'}:=\frac{tr(A^*D_{\Sigma}B)}{p^2}\ .$$
This first inner product is related to the loss function $l(.,.)$ through the identity $$l\left(\theta',\theta\right)=\|C(\theta')-C(\theta)\|^2_{\mathcal{H}}\ .$$ Besides, these two inner products clearly satisfy $\|C(\theta')\|_{\mathcal{H}}=\|D(\theta')\|_{\mathcal{H'}}$ for any $\theta'\in \Theta^+$. Gathering these new notations, we may upper bound (\ref{borne}) by
\begin{eqnarray}
\overline{\gamma}_{n,p}(\theta_{m,\rho_1}) & -& \overline{\gamma}_{n,p}(\widetilde{\theta}_{\rho_1})
    \leq    
  \|[I_{p^2}-D(\theta_{m,\rho_1})]^2-[I_{p^2}-D(\widetilde{\theta}_{\rho_1})]^2\|_{\mathcal{H'}}\times\nonumber\\
&&{\small  \sup_{{\tiny \begin{array}{c}
     \theta_1\in \Theta_m, \theta_2\in \Theta_{\widehat{m}},\\      \|[I_{p^2}-D(\theta_1)]^2-[I_{p^2}-D(\theta_2)]^2\|_{\mathcal{H'}}\leq 1 \end{array}}}
    \left<\left[I_{p^2}-D(\theta_1)\right]^2-\left[I_{p^2}-D(\theta_2)\right]^2,\left[\overline{\bf YY^*}-I_{p^2}\right]\right>_{\mathcal{H}'}\ .} \label{preuve3}
\end{eqnarray}

The first term in this product is easily bounded as these matrices are
diagonal. 
\begin{eqnarray}
\big\|\big[I_{p^2}-D(\theta_{m,\rho_1})\big]^2-\big[I_{p^2}-D(\widetilde{\theta}_{\rho_1})\big]^2\big\|_{\mathcal{H'}} =  tr\left[\left([I_{p^2}-D(\theta_{m,\rho_1})]^2-[I_{p^2}-D(\widetilde{\theta}_{\rho_1})]^2\right)^2\frac{D_{\Sigma}}{p^2}\right]^{\frac{1}{2}} \nonumber\\
  =  tr\left[\left[D(\theta_{m,\rho_1})-D(\widetilde{\theta}_{\rho_1})\right]^2\frac{D_{\Sigma}}{p^2} \left[2I_{p^2}-D(\theta_{m,\rho_1})-D(\widetilde{\theta}_{\rho_1})\right]^2 \right]^{1/2} \nonumber \\
 \leq  \varphi_{\text{max}}\left[ 2I_{p^2} -D(\theta_{m,\rho_1})-D(\widetilde{\theta}_{\rho_1}) \right] \|D(\theta_{m,\rho_1})-D(\widetilde{\theta}_{\rho_1})\|_{\mathcal{H'}}\ . \label{preuve4}
\end{eqnarray}
Since $\theta_{m,\rho_1}$ and $\widetilde{\theta}_{\rho_1}$ respectively belong to $\Theta_{m,\rho_1}^+$ and		 $\Theta_{\widehat{m},\rho_1}^+$, the largest eigenvalues of the matrices $I_{p^2}-C(\theta_{m,\rho_1})$ and $I_{p^2}-C(\widetilde{\theta}_{\rho_1})$ are smaller than $\rho_1$. Hence, we get  
\begin{eqnarray*}
 \varphi_{\text{max}}\left[ 2I_{p^2} -D(\theta_{m,\rho_1})-D(\widetilde{\theta}_{\rho_1})\right]  =  \varphi_{\text{max}}\left[ I_{p^2}-
C(\theta_{m,\rho_1})\right] + \varphi_{\text{max}}\left[ I_{p^2}-C(\widetilde{\theta}_{\rho_1})\right]   \leq  2\rho_1 \ .
\end{eqnarray*} 
Let us turn to the second term in (\ref{preuve3}). First, we
embed the set of matrices over which the supremum is taken in a ball of a
vector space. For any model $m'\in
\mathcal{M}$, let $U_{m'}$ be the space generated by the matrices $D(\theta')^2$ and $D(\theta')$ for  $\theta'\in \Theta_{m'}$. In the sequel, we note $d_{m'^2}$ the dimension of $U_{m'}$. The space $U_{m,m'}$ is defined as
the sum of $U_m$ and $U_{m'}$ whereas $d_{m^2,m'^2}$ stands for its dimension.
Finally, we note $\mathcal{B}^{\mathcal{H}'}_{m^2,m'^2}$ the unit ball of
$U_{m,m'}$ with respect to the inner product $<|>_{\mathcal{H}'}$. Gathering these notations, we get
$$\sup_{\small \begin{array}{c} R= \left[I-D(\theta_1)\right]^2-\left[I_{p^2}-D(\theta_2)\right]^2,\\\theta_1\in \Theta_m, \theta_2\in \Theta_{\widehat{m}}
    \text{ and }
 \|R\|_{\mathcal{H}'}\leq 1
    \end{array}}
\left<R,\overline{\bf YY^*}-I_{p^2}\right>_{\mathcal{H'}}
\leq\\
 \sup_{R\in \mathcal{B}^{\mathcal{H}'}_{m^2,\widehat{m}^2} }
\frac{1}{p^2}tr\left[RD_{\Sigma}\left(\overline{\bf YY^*}-I_{p^2}\right)\right]\ .$$
Applying the classical inequality $ab\leq \delta a^2+\delta^{-1}b^2/4$ and gathering
inequalities (\ref{preuve3}) and (\ref{preuve4}) yields
\begin{eqnarray}
\lefteqn{\overline{\gamma}_{n,p}(\theta_{m,\rho_1})-\overline{\gamma}_{n,p}(\widetilde{\theta}_{\rho_1}) \leq} & & \nonumber\\ &&\delta^{-1}
\|C(\theta_{m,\rho_1})-C(\widetilde{\theta}_{\rho_1})\|^2_{\mathcal{H}} + \rho_1^2\delta  \sup_{R\in
  \mathcal{B}^{\mathcal{H}'}_{m^2,\widehat{m}^2} }
\frac{1}{p^2}tr^2\left[RD_{\Sigma}\left(\overline{\bf YY^*}-I_{p^2}\right)\right]\ .\label{preuve5}
\end{eqnarray}
For any model $	m'\in \mathcal{M}$, we define the random variable $Z_{m'}$ as
\begin{eqnarray*}
Z_{m'} := \sup_{R\in \mathcal{B}^{\mathcal{H}'}_{m^2,m'^2}}
\frac{1}{p^2}tr\left[RD_{\Sigma}\left(\overline{\bf YY^*}-I_{p^2}\right)\right]\ .
\end{eqnarray*}
The variables $Z_{m'}$ turn out to be suprema of Gaussian chaos of order
2. In order to bound $Z_{\widehat{m}}$, we simultaneously control the deviations of $Z_{m'}$ for any model $m' \in \mathcal{M}$ thanks to the	 following lemma.

\begin{lemma}\label{concentrationZ}
For any positive numbers $\alpha$ and $\xi$ and any model $m'\in\mathcal{M}$, 
\begin{eqnarray*}
\mathbb{P}\left(Z_{m'} \geq \sqrt{\frac{2 \varphi_{\text{max}}(\Sigma)
 }{n}}\left\{\sqrt{1+\alpha/2}\sqrt{d_{m^2,m'^2}}+ \xi \right\}\right) \leq \hspace{6cm} \nonumber\\ \hspace{4cm}\exp\left\{-
 L_2\sqrt{d_{m'}}\left(\frac{\alpha}{\sqrt{1+\alpha/2}}\wedge
 \frac{\alpha^2}{1+\alpha/2}\right) -
 L_1\xi\left[\frac{\alpha}{\sqrt{1+\alpha/2}}\wedge \sqrt{n}\right] \right\}\ .
\end{eqnarray*}
\end{lemma}
This result is a consequence from a general concentration inequality for suprema
Gaussian chaos of order $2$ stated in Proposition \ref{lemmechaos}. Its proof is postponed to the technical appendix \cite{technical}. Let us fix the positive numbers $\alpha$ and $\xi$. Applying Lemma \ref{concentrationZ} to any model $m'\in\mathcal{M}$, the  event $\Omega'_\xi$ defined by 
$$\Omega'_{\xi}= \left\{ Z_{\widehat{m}}\leq \sqrt{\frac{2 \varphi_{\text{max}}(\Sigma) }{n}}\left[\sqrt{1+\alpha/2}\sqrt{d_{m^2,\widehat{m}^2}}+ \xi \right]\right\}$$
satisfies
$$\mathbb{P}(\Omega'^c_{\xi})\leq\exp\left\{-L_1\xi\left[\frac{\alpha}{\sqrt{1+\alpha/2}}\wedge \sqrt{n}\right] \right\} \sum_{m'\in\mathcal{M}}\exp\left\{-
 L_2\sqrt{d_{m'}}\left(\frac{\alpha}{\sqrt{1+\alpha/2}}\wedge
 \frac{\alpha^2}{1+\alpha/2}\right)\right\}\ . $$
From inequality (\ref{preuve5}), it follows that 
\begin{eqnarray}
\overline{\gamma}_{n,p}(\theta_{m,\rho_1})-\overline{\gamma}_{n,p}(\widetilde{\theta}_{\rho_1})
 & \leq &   \delta^{-1} \|C(\theta_{m,\rho_1})-C(\widetilde{\theta}_{\rho_1})\|^2_{\mathcal{H}} +\frac{2\delta \rho_1^2\varphi_{\text{max}}(\Sigma)}{np^2}   \left\{\sqrt{1+\alpha/2}\sqrt{d_{m^2,\widehat{m}^2}}+ \xi \right\}^2\ , \nonumber
\end{eqnarray}
conditionally to $\Omega'_{\xi}$. By triangle inequality, 
$$\|C(\theta_{m,\rho_1})-C(\widetilde{\theta}_{\rho_1})\|_{\mathcal{H}}\leq \|C(\theta_{m,\rho_1})-C(\theta)\|_{\mathcal{H}}+
\|C(\widetilde{\theta}_{\rho_1})-C(\theta)\|_{\mathcal{H}}\ .$$ 
We recall that  the loss function $l\left(\theta',\theta\right)$ equals $\|C(\theta')-C(\theta)\|^2_{\mathcal{H}}$. We apply twice the inequality  $(a+b)^2\leq (1+\beta)a^2+(1+\beta^{-1})b^2$. Setting the first $\beta$ to $\sqrt{\delta}-1$, it follows that
\begin{eqnarray*}
\overline{\gamma}_{n,p}(\theta_{m,\rho_1})-\overline{\gamma}_{n,p}(\widetilde{\theta}_{\rho_1})
 &\leq& \frac{1}{\sqrt{\delta}}l(\widetilde{\theta}_{\rho_1},\theta)+
\frac{\sqrt{\delta}}{\sqrt{\delta}-1} l(\theta_{m,\rho_1},\theta) \nonumber \\  &+&
\frac{2\delta \rho_1^2\varphi_{\text{max}}(\Sigma)}{np^2}\left[d_{m^2,\widehat{m}^2} (1+\beta)(1+\alpha/2)+\xi^2(1+\beta^{-1})\right]\ .
\end{eqnarray*}
By definition of $U_{m,\widehat{m}}$, its dimension $d_{m^2,\widehat{m}^2}$ is bounded by $d_{m^2}+
d_{\widehat{m}^2}$. Choosing $\beta=\delta-1$ yields 
\begin{eqnarray}
\overline{\gamma}_{n,p}(\theta_{m,\rho_1})-\overline{\gamma}_{n,p}(\widetilde{\theta}_{\rho_1})
 & \leq & \frac{1}{\sqrt{\delta}}l(\widetilde{\theta}_{\rho_1},\theta)+\frac{\sqrt{\delta}}{\sqrt{\delta}-1} l(\theta_{m,\rho_1},\theta) \label{preuve6}
 \\&& +
\frac{2\delta^2 \rho_1^2\varphi_{\text{max}}(\Sigma)}{np^2}\left[d_{m^2}
  (1+\alpha/2)+d_{\widehat{m}^2} (1+\alpha/2)\right] +
\frac{8\xi^2\varphi_{\text{max}}(\Sigma) \delta^2}{np^2(\delta-1)}\ . \nonumber 
\end{eqnarray}
To conclude, we need to compare the dimension $d_{m'^2}$ of the space $U_{m'}$ with $d_{m'}$.
\begin{lemma}\label{lemmedim}
For any model $m\in \mathcal{M}$, it holds that  
$$d_{m^2}\leq L d_{m}\ ,$$
where $L$ is a numerical constant between 4 and 5.48. 
\end{lemma}
The proof is postponed to the technical appendix~\cite{technical}. Defining the universal constant $K_0:=2L$, we derive from (\ref{preuve6}) that
\begin{eqnarray*}
\overline{\gamma}_{n,p}(\theta_{m,\rho_1})-\overline{\gamma}_{n,p}(\widetilde{\theta}_{\rho_1})
  &\leq & \frac{1}{\sqrt{\delta}}l(\widetilde{\theta}_{\rho_1},\theta)+\frac{\sqrt{\delta}}{\sqrt{\delta}-1} l(\theta_{m,\rho_1},\theta)  \\ & +&\frac{K_0\delta^2 \rho_1^2\varphi_{\text{max}}(\Sigma)}{np^2}\left[d_{m}
  (1+\alpha/2)+d_{\widehat{m}} (1+\alpha/2) +
\frac{\xi^2}{\delta-1}\right]\ ,
\end{eqnarray*}
with probability larger than $\mathbb{P}(\Omega'_{\xi})$. The isotropic case is analogous if we replace $d_m$ by $d_m^{\text{iso}}$.
\end{proof}\vspace*{0.5cm}

\subsection{Proofs of the minimax results}

Let us first prove a minimax lower bound on hypercubes $\mathcal{C}_m(\theta',r)$. We recall that these hypercubes are introduced in Definition \ref{definition_hypercube}. 

\begin{lemma}\label{minoration_hypercube}
Let $m$ be a model in $\mathcal{M}_1$ that satisfies $d_m\leq\sqrt{n}p$ and let $\theta'$ 
be a matrix in $\Theta_m\cap\mathcal{B}_1(0_p,1)$. Then, for any positive number $r$ such that $(1-\|\theta'\|_1-2rd_m)$ is positive, 
\begin{eqnarray*}
\inf_{\widehat{\theta}} \sup_{\theta \in
  \text{\emph{Co}}\left[\mathcal{C}_m(\theta',r)\right]}\mathbb{E}_{\theta}\left[l\left(\widehat{\theta},\theta\right)\right]
  \geq  L \sigma^2 \left(r\wedge \frac{1-\|\theta'\|_1}{\sqrt{np^2}}\right)^2 d_m\ ,
\end{eqnarray*}
where $\text{\emph{Co}}\left[\mathcal{C}_m(\theta',r)\right]$ denotes the convex hull of $C_m(\theta',r)$. Similarly, let $m$ be a model in $\mathcal{M}_1$ such $d_m^{\text{\emph{iso}}}\leq\sqrt{n}p$ and let $\theta'$ 
be a matrix in $\Theta_m^{\text{\emph{iso}}}\cap\mathcal{B}_1(0_p,1)$. Then, for any positive number $r$ such that $(1-\|\theta'\|_1-8rd_m^{\text{\emph{iso}}})$ is positive, 
\begin{eqnarray*}
\inf_{\widehat{\theta}} \sup_{\theta \in
  \text{\emph{Co}}\left[\mathcal{C}^{\text{\emph{iso}}}_m(\theta',r)\right]}\mathbb{E}_{\theta}\left[l\left(\widehat{\theta},\theta\right)\right]
  \geq  L \sigma^2 \left(r\wedge \frac{1-\|\theta'\|_1}{\sqrt{np^2}}\right)^2 d^{\text{\emph{iso}}}_m\ .
\end{eqnarray*}
\end{lemma}

\begin{proof}[Proof of Proposition \ref{variance_minimax_global}]
The first result derives from Lemma \ref{minoration_hypercube} applied to the hypercube $\mathcal{C}_m(0_p,(np^2)^{-1/2})$. We prove the second result using the same lemma with  $\mathcal{C}_m[\theta',(1-\|\theta\|_1)/(\sqrt{n}p)]$.
\end{proof}\vspace{0.5cm}

\begin{proof}[Proof of Lemma \ref{minoration_hypercube}]
This lower bound is based on an application of Fano's approach. See \cite{yu} for a review of
this method and comparisons with Le Cam's and Assouad's Lemma. The proof follows three main steps:
First, we upper bound the Kullback-Leibler entropy between distributions
corresponding to $\theta_1$ and $\theta_2$ in the hypercube. Second, we find a set of points in the hypercube well
separated with respect to the Hamming distance. Finally, we conclude by applying Birg\'e's version
of Fano's lemma. More details can be found in the technical appendix ~\cite{technical}. 

\end{proof}
\vspace{0.5cm}

\begin{proof}[Proof of Proposition \ref{minoration_ellipsoid}]

First, observe that the set $\mathcal{E}(a)\cap\mathcal{B}_1(0_p,1/2)$  is included in $\mathcal{E}(a)\cap\mathcal{B}_1(0_p,1)\cap\mathcal{U}(2)$. We then derive minimax lower bounds on $\mathcal{E}(a)\cap\mathcal{B}_1(0_p,1/2)$ from the lower bounds on hypercubes.

 Let $m_i$ be a model in $\mathcal{M}_1$ such that $d_m$ is smaller than $\sqrt{n}p$.
Let us look for positive numbers $r$ such that the hypercube $\left[\mathcal{C}_{m_i}(0_p,r)\right]$ is included in
the set $\mathcal{E}(a)\cap\mathcal{B}_1(0_p,1/2)$. 
\begin{lemma}\label{majoration_variance}
Let $m$ be a model in $\mathcal{M}_1$ and $r$ be a positive number smaller than $1/(4d_m)$. For any $\theta\in \text{\emph{Co}}\left[\mathcal{C}_{m}(0_p,r)\right]$, 
\begin{eqnarray}
\text{\emph{var}}_{\theta}\left(X{\scriptstyle[0,0]}\right) \leq \sigma^2\left(1 + 16d_mr^2\right)\ . \nonumber
\end{eqnarray}
\end{lemma}
The proof is postponed to the technical appendix \cite{technical}. If we choose 
$$r\leq \frac{a_i}{16\sigma\sqrt{d_{m_i}}}\ ,$$
then $2rd_{m_i}$ is smaller than $1/8$ by assumption $(\mathbb{H}_a)$. Applying Lemma \ref{majoration_variance}, we then derive that  $\var_{\theta}\left(X{\scriptstyle[0,0]}\right) \leq \sigma^2 + a_i^2$. Hence,  we get the upper bound ~\\ $\sum_{j=1}^i\left[
\var\left(X{\scriptstyle[0,0]}|X_{m_{j-1}}\right)-\var\left(X{\scriptstyle[0,0]}|X_{m_j}\right)\right]\leq
a_i^2$ and it follows that  
\begin{eqnarray*}
 \sum_{j=1}^{\text{Card}(\mathcal{M}_1)}\frac{
\var\left(X{\scriptstyle[0,0]}|X_{m_{k-1}}\right)-\var\left(X{\scriptstyle[0,0]}|X_{m_j}\right)}{a_j^2}\leq 1\ ,
\end{eqnarray*}
 since the sequence $(a_j)_{1\leq j \leq \text{Card}(\mathcal{M}_1)}$ is non increasing.
Consequently, $\text{Co}\left[\mathcal{C}_m(0_p,r)\right]$ is a subset of $\mathcal{E}(a)\cap\mathcal{B}_1(0_p,1/2)$. By Lemma \ref{minoration_hypercube}, we get
\begin{eqnarray}
\inf_{\widehat{\theta}}\sup_{\theta\in \mathcal{E}(a)\cap\mathcal{B}_1(0_p,1/2)}
\mathbb{E}_{\theta}\left[l\left(\widehat{\theta},\theta\right)\right]& \geq &
L \sigma^2\left(\frac{a_i^2}{16\sigma^2}\wedge \frac{d_{m_i}}{np^2}\right)\nonumber \\
& \geq & L \left(a_i^2\wedge \frac{\sigma^2d_{m_i}}{np^2}\right) \ .\label{majoration_quotient_dimension2}
\end{eqnarray}
Considering all models $m\in\mathcal{M}_1$ such that $d_m\leq \sqrt{n}p$ yields
\begin{eqnarray}
 \inf_{\widehat{\theta}}\sup_{\theta\in \mathcal{E}(a)\cap\mathcal{B}_1(0_p,1/2)}
\mathbb{E}_{\theta}\left[l\left(\widehat{\theta},\theta\right)\right]& \geq & L\sup_{i\leq \text{Card}(\mathcal{M}_1), \ d_{m_i}\leq \sqrt{n}p} \left(a_i^2\wedge \frac{\sigma^2d_{m_i}}{np^2}\right) \ .\label{minoration_el}
\end{eqnarray}
If the maximal dimension $d_{m_{\text{Card}(\mathcal{M}_1)}}$ is smaller than $\sqrt{n}p$, the proof is finished. In the opposite case, we need to show that the supremum (\ref{minoration_ellipsoid_equation}) over all models $m\in \mathcal{M}_1$ is achieved at some model $m$ of dimension less than $\sqrt{n}p$.
\begin{lemma}\label{majoration_quotient_dimension}
For any integer $1\leq i\leq \text{Card}(\mathcal{M}_1)-1$, the ratio $d_{m_{i+1}}/d_{m_i}$ is less than  $2$.
\end{lemma}
The proof of Lemma \ref{majoration_quotient_dimension} is postponed to the technical appendix \cite{technical}.
Let $i'$ be the largest integer such that $d_{m_{i'}}\leq \sqrt{n}p$. Since
$i'$ is smaller than $\text{Card}(\mathcal{M}_1)$, we know from Lemma \ref{majoration_quotient_dimension} that $\sqrt{n}p/2\leq
d_{m_{i'}}\leq \sqrt{n}p$. By assumption $(\mathbb{H}_a)$, $a_{i'}^2$ is smaller than $\sigma^2/d_{m_{i'}}$.
Gathering these bounds yields 	
$$a_{i'}^2\leq \frac{\sigma^2}{d_{m_{i'}}}\leq \frac{4d_{m_{i'}}\sigma^2}{np^2}\ .$$
Since the sequence $(a_i)_{1\leq i\leq \text{Card}(\mathcal{M}_1)}$ is non increasing, the supremum  (\ref{minoration_ellipsoid_equation}) over all models in $\mathcal{M}_1$ is either achieved for some $i\leq i'$ or is smaller than $4(a_{i'}^2\wedge \sigma^2d_{m_{i'}}/(np^2))$. 
\end{proof} \vspace{0.5cm}

\begin{proof}[Proof of Corollary \ref{adaptation_sparsite}]
 Observe that $\text{Co}[\mathcal{C}_m(0_p,1/(4d_m)]$ is included in $\Theta_{m}\cap\mathcal{B}_1(0_p,1/2)$. This last set is itself included in  $\Theta_{m,\rho_1}^+\cap\mathcal{U}(\rho_2)$. Applying Lemma \ref{minoration_hypercube}, we get the following minimax lower bound
\begin{eqnarray*}		
 \inf_{\widehat{\theta}}\sup_{\theta\in\Theta^+_{m,\rho_1}\cap\mathcal{U}(\rho_2)}\mathbb{E}\left[l\left(\widehat{\theta},\theta\right)\right]\geq L\sigma^2\frac{d_m}{np^2}\ ,
\end{eqnarray*}
since the dimension $d_m$ is smaller than $np^2$.
Applying Theorem \ref{mainthrm}, we derive that
\begin{eqnarray*}
 \sup_{\theta\in\Theta^+_{m,\rho_1}\cap\mathcal{U}(\rho_2)}\mathbb{E}\left[l\left(\widetilde{\theta}_{\rho_1},\theta\right)\right] & \leq & L(K)\sigma^2\rho_1^2\rho_2\frac{d_m}{np^2}+L_2(K)\frac{\rho_1^2}{np^2}\sup_{\theta\in\Theta^+_{m,\rho_1}\cap\mathcal{U}(\rho_2)}\varphi_{\text{max}}(\Sigma)\\ 
&\leq & L(K,\rho_1,\rho_2)\sigma^2\frac{d_m}{np^2} .
\end{eqnarray*}
We conclude by combining the two different bounds.
\end{proof}

\begin{proof}[Proof of Proposition \ref{prte_adaptatif}]
This result derives from the upper bound of the risk of $\widetilde{\theta}_{\rho_1}$ stated in Theorem \ref{mainthrm} and the minimax lower bound stated in Proposition \ref{minoration_ellipsoid}. For details, we refer to the technical appendix~\cite{technical}.

\end{proof}

\subsection{Proofs of the asymptotic risk bounds}\label{variance_proofs}

\begin{proof}[Proof of Proposition \ref{variance_vrai}]
This result is closely related to Proposition 4.11 in \cite{guyon95}. In fact, we extend
his proof to stationary fields on a torus. In the sequel, we shall only consider non-isotropic GMRFs, the isotropic case being similar. Let us fix a model $m$ in the collection $\mathcal{M}_1$ and let us assume $(\mathbb{H}_1)$.

We define the $d_m\times p^2$ matrix $\chi_m^v$ as
$$\left(\chi_m^v\right)^* := \left(\left[C(\Psi_{i_k,j_k})X^v\right],\ k=1,\ldots, d_m\right)\ . $$
For any $(i,j)\in \{1,\ldots,p\}^2$, the $(i-1)p+j$-th row of $\chi_m^v$ corresponds to the list of covariates used when performing the regression of $X{\scriptstyle[i,j]}$ with respect to its neighbours  in the model $m$. Contrary to the previous proofs, we need to express the $n\times p^2$ matrix ${\bf X^v}$ in terms of a vector. This is why we define the vector ${\bf X^V}$ of size $np^2$ as
\begin{eqnarray*}
{\bf X^V}{\scriptstyle[p^2(j-1)+p(i_1-1)+i_2]} := {\bf X}^j{\scriptstyle[i_1,i_2]}\ ,
\end{eqnarray*}
for any $(i_1,i_2)\in \{1,\ldots,p\}^2$ and any $j\leq n$. Similarly, let $\boldsymbol{\chi}_m^V$ be the $d_m\times np^2$ matrix defined as 
\begin{eqnarray*}
\boldsymbol{\chi}_m^{\bf V}{\scriptstyle[k,p^2(j-1)+p(i_1-1)+i_2]} := \boldsymbol{\chi}_m^{j}{\scriptstyle[p(i_1-1)+i_2]}\ , 
\end{eqnarray*}
for any $(i_1,i_2)\in \{1,\ldots,p\}^2$ and any $j\leq n$.

We are not able to work out directly the asymptotic risk of $\widehat{\theta}_{m,\rho_1}$. This is why we introduce a new estimator $\check{\theta}_m$ whose asymptotic distribution is easier to derive. Afterwards, we shall prove that $\check{\theta}_m$ and $\widehat{\theta}_{m,\rho_1}$ have the same asymptotic distribution. Let us respectively define the estimators $\check{a}_m$ in $\mathbb{R}^{d_m}$ and
$\check{\theta}_m$ as
\begin{eqnarray}\label{definition_a}
\check{a}_m&  := & \left(\left(\boldsymbol{\chi}_m^{{\bf V}}\right)^*
\boldsymbol{\chi}_m^{{\bf V}}\right)^{-1}\boldsymbol{\chi}_m^{\bf V}{\bf X^V}\\
\check{\theta}_m & := & \sum_{k=1}^{d_m} \check{a}_m{\scriptstyle[k]}\Psi_{i_k,j_k} \nonumber\ ,
\end{eqnarray}
where we recall that $(\Psi_{i_1,j_1},\ldots ,\Psi_{i_{d_m},j_{d_m}})$ is a basis of $\Theta_m$. Obviously, $\check{\theta}_m$ is a Conditional least squares estimator since it minimizes the expression (\ref{definition_gamman}) of $\gamma_{n,p}(.)$ over the whole space $\Theta_m$.
Consequently, $\check{\theta}_m$ coincides with $\widehat{\theta}_{m,\rho_1}$ if $\check{\theta}_m$ belongs to $\Theta_{m,\rho_1}^+$.

For the second result, we assume that Assumption $(\mathbb{H}_2)$ holds. Applying Corollary \ref{definition_thetam}, we know that for any $(k,l)\in \Lambda$, $X{\scriptstyle[k,l]}$ decomposes as
\begin{eqnarray}\label{ecriture_conditionnelle}
X{\scriptstyle[k,l]} = \sum_{(i,j)\in  m }\theta_{m,\rho_1}{\scriptstyle[i,j]}X{\scriptstyle[k+i,l+j]} + \epsilon_m{\scriptstyle[k,l]}\ ,
\end{eqnarray}
where $\epsilon_m{\scriptstyle[k,l]}$ is independent from $\left\{X{\scriptstyle[k+i,l+j]},\  (i,j)\in m \right\}$. For the first result, the same decomposition holds since $\theta$ is assumed to belong to $\Theta_{m,\rho_1}^+$ and $\theta_{m,\rho_1}$ therefore equals $\theta$.

Let $a_m\in\mathbb{R}^{d_m}$ be the unique vector such that
$\theta_{m,\rho_1} = \sum_{k=1}^{d_m}a_m{\scriptstyle[k]}\Psi_{i_k,j_k}$. Then, the previous decomposition  becomes 
\begin{eqnarray*}
X^v = a_m^*\chi_m^v +\epsilon_m^v\ .
\end{eqnarray*}
Gathering this last identity with (\ref{definition_a}) yields 
\begin{eqnarray*}
\check{a}_m - a_m  = \left(\frac{1}{np^2}(\boldsymbol{\chi}_m^{\bf V})^*
\boldsymbol{\chi}_m^{\bf V}\right)^{-1}\left(\frac{1}{np^2}\boldsymbol{\chi}_m^{\bf V}\boldsymbol{\epsilon}_m^{\bf V}\right)\ ,
\end{eqnarray*}
where the vector $\boldsymbol{\epsilon}_m^{\bf V}$ of size $np^2$ corresponds to the $n$ observations of the vector $\epsilon_m^v$.
When $n$ goes to the infinity, $1/(np^2) (\boldsymbol{\chi}_m^{\bf V})^*
\boldsymbol{\chi}_m^{\bf V}$ converges almost surely to the covariance matrix $V$ by the law of large numbers. By definition, the variable $\epsilon_m{\scriptstyle[i,j]}$ is independent from the $(i-1)p+j$th row of $\chi^v_m{\scriptstyle[i,j]}$. It follows that
$\mathbb{E}_{\theta}( \boldsymbol{\chi}_m^{\bf V}\boldsymbol{\epsilon}^{\bf V}
)=0$. Applying again the law of large numbers we conclude that
$\check{a}_m$ converges almost surely towards $a_m$ and that
$\check{\theta}_m$ converges almost surely towards $\theta_{m,\rho_1}$. Besides, the central limit theorem states that the random vector
$1/(\sqrt{n}p)\boldsymbol{\chi}_m^{\bf V}\boldsymbol{\epsilon}^{\bf V}$
 converges in distribution towards a zero mean Gaussian vector
whose covariance matrix equals $1/p^2\var_{\theta}\left(
\chi_m^v\epsilon_m^{v}\right)$. By decomposition (\ref{ecriture_conditionnelle}), $\epsilon_m^v=(I-C(\theta_{m,\rho_1}))X^v$ while the $k$-th row of $\chi_m^v$ equals $\left[C(\Psi_{i_k,j_k})X^v\right]^*$. Thus, for any $1\leq k,l\leq d_m$,  
\begin{eqnarray*}
 \frac{1}{p^2}\var_{\theta}\left(
\chi_m^v\epsilon_m^{v}\right){\scriptstyle[k,l]} = \frac{1}{p^2}\cov_{\theta}\left[(X^v)^*C(\Psi_{i_k,j_k})\left[I-C(\theta_{m,\rho_1})\right]X^v,(X^v)^*C(\Psi_{i_l,j_l})\left[I-C(\theta_{m,\rho_1})\right]X^v\right]\ .
\end{eqnarray*}
As the covariance matrix of $X^v$ is $\sigma^2\left(I-C(\theta)\right)^{-1}$, we obtain by standard Gaussian properties
\begin{eqnarray*}
\frac{1}{p^2}\var_{\theta}\left(
\chi_m^v\epsilon_m^{v}\right){\scriptstyle[k,l]}=\hspace{10cm} \nonumber \\ \hspace{1cm} \frac{2\sigma^4}{p^2}\cov_{\theta}\left[\left[I-C(\theta)\right]^{-1}C(\Psi_{i_k,j_k})\left[I-C(\theta_{m,\rho_1})\right]\left[I-C(\theta)\right]^{-1}C(\Psi_{i_l,j_l})\left[I-C(\theta_{m,\rho_1})\right]\right]\ . 
\end{eqnarray*}
By Lemma \ref{codiagonalisation}, all these matrices are diagonalizable in the same basis and therefore commute with each other. We conclude that $\frac{1}{p^2}\var_{\theta}\left(
	\chi_m^v\epsilon_m^{v}\right)= 2\sigma^4W$ and 
\begin{eqnarray*}
\sqrt{n}p\left(\check{a}_m-a_m\right)\rightarrow \mathcal{N}\left(0,V^{-1}WV^{-1}\right)\ .
\end{eqnarray*}

As $\widehat{\theta}_{m,\rho_1}$ belongs to $\Theta^{+}_{m,\rho_1}$, there exists a unique vector
$\widehat{a}_m\in \mathbb{R}^{d_m}$ such that
$\widehat{\theta}_{m,\rho_1} = \sum_{k= 1}^{d_m}\widehat{a}_m{\scriptstyle[k]}\Psi_{i_k,j_k}$.
 The matrix $\theta_{m,\rho_1}$ belongs to the open set $\Theta^{+}_{m,\rho_1}$ for the two cases of the propositions. Indeed, $\theta_{m,\rho_1}$ equals $\theta$ in the first situation. In the second situation, this is due to the fact that $\theta$ satisfies $(\mathbb{H}_2)$ and to Lemma \ref{lemme_projection}.

Since $\check{\theta}_m$ converges
almost surely to $\theta_{m,\rho_1}$, the matrix $\check{\theta}_m$ belongs to $m$ with probability going to one when $n$ goes to infinity. If follows that the estimators  
$\check{a}_m$ and $\widehat{a}_m$ coincide with probability going to one. By Slutsky's Lemma, we obtain that  	 
\begin{eqnarray*}
\sqrt{n}p\left(\widehat{a}_m-a_m\right) \rightarrow \mathcal{N}\left(0,V^{-1}WV^{-1}\right)\ .
\end{eqnarray*}
Let us express the risk of $\widehat{\theta}_{m,\rho_1}$ with respect to the distribution of $\widehat{a}_m$.
\begin{eqnarray*}
l\left(\widehat{\theta}_{m,\rho_1},\theta_{m,\rho_1}\right) = \mathbb{E}_{\theta}\bigg[\sum_{k=1}^{d_m}
\left(\widehat{a}_m{\scriptstyle[k]}-a_m{\scriptstyle[k]}\right)tr\left(\Psi_{i_k,j_k}X\right)\bigg]^2  = tr\left[V\left(\widehat{a}_m-a_m\right)^*\left(\widehat{a}_m-a_m\right)\right]\ .
\end{eqnarray*}
By Portmanteau's Lemma, $np^2l(\widehat{\theta}_{m,\rho_1},\theta_{m,\rho_1})$ converges in distribution towards a random variable whose expectation is $tr\left(WV^{-1}\right)$. In order to conclude, it remains to prove that the sequence $[np^2l(\widehat{\theta}_{m,\rho_1},\theta)	]_{n\geq 1}$ is asymptotically uniformly integrable. 

Let us consider a model selection procedure with the collection $\mathcal{M}=\{m\}$ and a penalty term satisfying the assumptions of Theorem \ref{mainthrm}. Arguing as in the proof of this theorem,  we derive from identity (\ref{equiintegrable}) the following property. For any $\xi>0$, with probability larger than $1-L_1\exp\left[-L_2\xi\right]$, 
\begin{eqnarray*}
 np^2l\left(\widehat{\theta}_{m,\rho_1},\theta_{m,\rho_1}\right)\leq L_3 d_m\varphi_{\text{max}}(\Sigma) + L_4\xi^2\varphi_{\text{max}}(\Sigma) \ .
\end{eqnarray*}
This clearly implies that the sequence $[np^2l(\widehat{\theta}_{m,\rho_1},\theta_{m,\rho_1})]_{n\geq 1}$ is asymptotically uniformly integrable and the first part of the result follows.\\

For the first result of the proposition, we have stated that $\theta$ equals $\Theta_m$. As a consequence,  
 $$\lim_{n\rightarrow +\infty}\mathbb{E}_{\theta}\left[l\left(\widehat{\theta}_{m,\rho_1},\theta\right)\right]=2\sigma^4tr\left[WV^{-1}\right]\ .$$
 Besides, the term $W{\scriptstyle[k,l]}$ here equals $tr\left[C(\Psi_{i_k,j_k})C(\Psi_{i_l,j_l})\right]$. This last quantity is zero if $k\neq l$ and equals $\|C(\Psi_{i_k,j_k})\|_F^2$ if $k=l$.
\end{proof}\vspace{0.5cm}

\begin{proof}[Proof of Proposition \ref{oracle_asymptotique}]
 As $\theta$ belongs to $\Theta^+\cap \mathcal{B}_1(0_p,\eta)$, the largest eigenvalue of $\Sigma$ is smaller than $\sigma^2/(1-\eta)$. Applying Theorem \ref{mainthrm}, we get
\begin{eqnarray*}
 \mathbb{E}_{\theta}\left[l\left(\widetilde{\theta}_{\rho_1},\theta\right)\right]& \leq & L(K) \inf_{m\in\mathcal{M}}\left[l(\theta_{m,\rho_1},\theta)+ K\frac{\sigma^2}{np^2(1-\eta)}\right]\\
& \leq & L(K,\eta) \inf_{m\in\mathcal{M}}\left[l(\theta_{m,\rho_1},\theta)+ K\frac{\sigma^2}{np^2}(1-\eta)^3\right] .
\end{eqnarray*}
Gathering this bound with the result of Corollary \ref{minoration_variance_asymptotique} enable us to conclude.
\end{proof}

\appendix

\section{}

\begin{lemma}\label{codiagonalisation}
There exists an orthogonal matrix $P$ which simultaneously diagonalizes every $p^2\times
p^2$ symmetric block circulant matrices with $p\times p$ blocks. Conversely, if  $\theta$
is a square matrix of size $p$ which satisfies (\ref{symmetric_matrice}), then 
the matrix $D(\theta)= PC(\theta)P^*$ is diagonal and satisfies
\begin{eqnarray}\label{forme_valeurs_propres}
D(\theta) {\scriptstyle[(i-1)p+j,(i-1)p+j]} = \sum_{k=1}^{p}\sum_{l=1}^{p}\theta{\scriptstyle[k,l]}\cos\left(2\pi(ki/p+lj/p)\right)
\end{eqnarray}
for any $1\leq i,j\leq p$. 
\end{lemma}
It is proved as in \cite{rue} Sect.2.6.2 to the price of a slight modification in order to take into account the fact that $P$ has is orthogonal and not unitary. The difference comes from the fact that contrary to Rue and Held we also assume that $C(\theta)$ is symmetric.

 This lemma states that all symmetric block circulant matrices are simultaneously diagonalizable. Moreover, Expression (\ref{forme_valeurs_propres}) explicitly provides the eigenvalues of the $C(\theta)$ as the two-dimensional discrete Fourier transform of the $p\times p$ matrix $\theta$.

\section*{Acknowledgements}
I am grateful to Pascal Massart for many fruitful discussions. I also thank the referees and the associate editor for their suggestions that led to an improvement of the manuscript.

\addcontentsline{toc}{section}{References}

\bibliographystyle{alpha}

\bibliography{spatial}

\end{document}